\documentclass[12pt]{article}

\usepackage{graphicx}
\usepackage{latexsym,amssymb}
\usepackage{amsthm}
\usepackage{indentfirst}
\usepackage{amsmath}
\usepackage{color}
\usepackage{xcolor}
\usepackage{fourier}
\usepackage[colorlinks=true,backref=page]{hyperref}
\usepackage[all]{xy}

\usepackage[a4paper,text={160true mm,230true mm},centering,top=35true mm,head=5true mm,headsep=2.5true mm,foot=8.5true mm]{geometry}

\newtheorem{theoremalph}{Theorem}

\newtheorem*{Main Theorem}{Main Theorem}
\newtheorem{Coro}[theoremalph]{Corollary}
\newtheorem{Theorem}{Theorem}[section]
\newtheorem*{Theorem A}{Theorem A}
\newtheorem*{Theorem A'}{Theorem A'}
\newtheorem*{Theorem B'}{Theorem B'}

\newtheorem*{Conjecture}{Conjecture}

\newtheorem{Proposition}[Theorem]{Proposition}
\newtheorem{Lemma}[Theorem]{Lemma}

\newtheorem{Remark}[Theorem]{Remark}

\newtheorem{Remark-numbered}[Theorem]{Remark}
\newtheorem{Corollary}[Theorem]{Corollary}

\newtheorem*{Claim}{Claim}
\newtheorem{Claim-numbered}[Theorem]{Claim}

 \def\NN{{\mathbb N}} 

 \def\RR{{\mathbb R}}

\def\cC{{\cal C}}    
   \def\cP{{\cal P}} 
   \def\cQ{{\cal Q}}

\newcommand{\loc}{\operatorname{loc}}

\def\dim{\operatorname{dim}}

\def\diam{\operatorname{Diam}}

\def\Jac{\operatorname{Jac}}

\begin{document}

\title{Upper semi-continuous properties for $\mathcal C^{1,\alpha}$ diffeomorphisms}

\author{Chiyi Luo and Dawei Yang\footnote{
D. Yang  was partially supported by National Key R\&D Program of China (2022YFA1005801), NSFC 12171348 \& NSFC 12325106, ZXL2024386 and Jiangsu Specially Appointed Professorship.
C. Luo was partially supported by NSFC (12501244) and Early-Career Young Scientists and Technologists Project of Jiangxi Province. 
}}

\date{}
\maketitle

\begin{abstract} 
We establish a uniform approximation of metric entropy by partition entropy using uniform partitions for $\cC^{1,\alpha}$ three-dimensional diffeomorphisms. 
This gives several consequences for diffeomorphisms on a compact manifold $M$ with ${\rm dim} M\leq 3$.

First, if an invariant measure $\mu$ is a continuity point of the sum of its positive Lyapunov exponents, then $\mu$ is an upper semi-continuity point of the entropy map. 
Second, it provides a slight improvement of the necessary condition for strong positive recurrence of surface and three-dimensional diffeomorphisms. 
Third, it yields continuity of dimensions for measures of maximal entropy.
\end{abstract}

\tableofcontents

\section{Introduction}\label{SEC:1}
The space of invariant measures may be very complicated for chaotic systems.  
The metric entropy of an invariant measure, which was first studied by Kolmogorov and Sinai, is a fundamental concept in ergodic theory.  
Denote by $h_\mu(f)$ the metric entropy of $\mu$ for a map $f$. 
The dependence of this quantity with respect to the invariant measures and the maps received people's great interest.
Because it is closely related to the existence of equilibrium states and measures of maximal entropy, which serve as fundamental objects for describing the statistical and dynamical behavior of the system.

In general, metric entropy is not continuous with respect to the measures.
Even for very nice system, it will not be lower-semi continuous with respect to the measures. 
For example, in a hyperbolic basic set with positive entropy, the measure of maximal entropy can be approximated by measures supported on periodic orbits. 
This leads to a failure of lower semi-continuity.

The upper semi-continuity of metric entropy holds in the uniformly hyperbolic setting, ensuring the existence of measures of maximal entropy. 
Another well-known class of systems with upper semi-continuous metric entropy is $\mathcal C^{\infty}$ diffeomorphisms.
Inspired by Yomdin's work \cite{Yom87} on Shub's entropy conjecture, Newhouse \cite{New89} proved for $\mathcal C^\infty$ differentiable maps, the upper semi-continuity of metric entropy holds with respect to invariant measures.
However, for $\mathcal C^r$ diffeomorphisms on a compact manifold with finite positive $r$, the upper semi-continuity of the metric entropy may fail, see, for instance, the counterexample in \cite{Buz14,Mis73}.

Let $f$ be a diffeomorphism of a compact manifold $M$, and let $\{\mathcal P_m\}_{m\ge 1}$ be a sequence of finite partitions such that $\diam(\cP_m) \to 0$.
Then, by the Kolmogorov–Sinai theorem, for every $f$-invariant measure $\mu$, we have $h_{\mu}(f,\cP_m) \to h_{\mu}(f)$.
Note that for a finite partition $\cP$, the map $\mu \mapsto h_{\mu}(f,\cP)$ is upper semi-continuous at every measure $\mu$ satisfying $\mu(\partial \cP)=0$.
Therefore, it is useful to obtain a uniform estimate for $\sup_{\mu} \big(h_{\mu}(f)-h_{\mu}(f,\cP)\big) $.
For instance, for expansive diffeomorphisms, there exists an $\varepsilon>0$ such that for every partition $\cP$ with $\diam(\cP)<\varepsilon$, one has $h_{\mu}(f)=h_{\mu}(f,\cP)$.
Then, the upper semi-continuity of the entropy map follows immediately.

Assume that $f:M\to M$ is a $\cC^{r,\alpha}$ diffeomorphism of a compact manifold $M$.
A deep result of Newhouse \cite{New89}, based on Yomdin's estimate \cite{Yom87}, shows that there exists $C_r>0$ such that for every $q\in \NN$, there is $\varepsilon_q>0$ such that 
$$h_\mu(f) \le h_\mu(f,\mathcal Q) + \frac{\dim M \cdot \log \|Df^q\|_{0}}{q(r+\alpha)} + \frac{\log(C_{r})}{q},$$
for every $\cQ$ with $\diam (\cQ)<\varepsilon_q$, where $\|Df^q\|_{0}:=\sup_{x\in M}\|D_xf^q\|$.
In particular, when $f:M\to M$ is $C^{\infty}$, this implies the upper semi-continuity of the entropy map $\mu \to h_{\mu}(f)$.
However, in the $C^r$ setting, Newhouse’s estimate only provides an upper bound for the possible failure of upper semi-continuity of the entropy map.

In this paper, we consider a $\cC^{1,\alpha}$ diffeomorphism $f:M\to M$ with $\alpha\in (0,1]$. 
We obtain a new upper bound in the case $\dim M\leq 3$ and demonstrate several applications of this estimate.

\subsection{Statement of the main estimates}
Suppose that $M$ is a compact Riemannian manifold without boundary and let $f:M\to M$ be a diffeomorphism. 
By the Oseledets theorem \cite{Ose68}, there exists an invariant set $\Gamma\subset M$ with total probability, 
i.e., $\mu(\Gamma)=1$ for any invariant measure $\mu$, such that for any $x\in\Gamma$, there are
\begin{itemize}
\item  $\lambda_1(x,f)>\cdots>\lambda_{s(x)}(x,f)$, which are $f$-invariant and measurable functions of $x$;
\item a $Df$-invariant measurable splitting $T_xM=E_1(x,f)\oplus E_2(x,f)\oplus \cdots\oplus E_{s(x)}(x,f)$
\end{itemize}
such that for any $v\in E_i(x)\setminus\{0\}$ with $1\le i\le s(x)$, one has
$$\lim_{n\to\pm\infty}\frac{1}{n}\log\|D_xf^n(v)\|=\lambda_i(x,f).$$
These numbers $\{\lambda_i(x,f)\}_{i=1}^{s(x)}$ are called the \emph{Lyapunov exponents} of $ x\in \Gamma$.

Given $x\in\Gamma$, we denote 
$$\max \big\{0,~\lambda_{1}(x,f) \big\}=\lambda^{+}(x,f),~~\min\big\{0,~\lambda_{s(x)}(x,f)\big\}=\lambda^{-}(x,f).$$
For each $f$-invariant measure $\mu$, define $\lambda^{\pm}(\mu,f)=\int \lambda^{\pm}(x,f)~{\rm d}\mu(x)$.
The sum of positive Lyapunov exponents is an important quantity used to describe the complexity of the dynamics. 
We define it as
$$\lambda_\Sigma^+(x,f):=\sum_{i=1}^{s(x)}\dim E_i(x) \cdot \max\{0,\lambda_i(x,f)\};$$
For an invariant measure $\nu$ of $f$, define the \emph{the sum of positive Lyapunov exponents} of $\nu$ by
$\lambda_\Sigma^{+}(\nu,f)=\int \lambda_\Sigma^+(x,f){\rm d}\nu(x)$.
Note that $\lambda^+_{\Sigma}(\mu, f)$ is upper semi-continuous with respect to $\mu$ and $f$. 
A detailed proof can be found in Proposition \ref{Pro:usc-sum-positive}.

We say that an invariant measure $\mu$ of $f$ has exactly one positive Lyapunov exponent, if $\mu$-almost every $x$ satisfies $\lambda_1(x,f)>0\geq \lambda_2(x,f)$ and $\dim E_1(x,f)=1$; and $\mu$ has exactly one negative Lyapunov exponent, if $\mu$-almost every $x$ satisfies $\lambda_{s(x)-1}(x,f)\geq 0> \lambda_{s(x)}(x,f)$ and $\dim E_{s(x)}(x,f)=1$.
The main theorem is stated as follows:
\begin{theoremalph}\label{Thm:A-entropy-bound-partition}
	Given $\alpha\in (0,1]$, there exists a constant $C_{\alpha}$ with the following property.
	For each $\Upsilon>0$ and $q\in\mathbb N$, there exists $\varepsilon_q=\varepsilon_{\Upsilon,q}>0$ such that 
	\begin{itemize}
		\item  for every $\mathcal{C}^{1,\alpha}$ diffeomorphism $f:~M\to M$ satisfying $\|f\|_{\mathcal C^{1,\alpha}}\leq \Upsilon$;
		\item  for every ergodic measure $\mu$ of $f$ with exactly one positive Lyapunov exponent;
		\item  for every finite partition $\mathcal Q$ with ${\rm Diam}(\mathcal Q)<\varepsilon_q$;
	\end{itemize}
	one has 
	$$h_\mu(f)\le h_\mu(f,\mathcal Q)+\frac{1}{\alpha}\big[\frac{1}{q}\int\log\|D_xf^q\|{\rm d}\mu(x)-\lambda^{+}(\mu,f)+\frac{1}{q}\big]+\frac{\log(2q \Upsilon C_{\alpha})}{q}.$$
\end{theoremalph}

In the case $\dim M \leq 3$, we obtain a uniform estimate for every invariant measure $\mu$. 
For $1\leq k\leq \dim M$, let $\wedge^k D_xf$ denote the induced action of $D_xf$ on the $k$-th exterior power.

\begin{theoremalph}\label{Thm:B-Invar}
	Assume that $\dim M\le 3$. 
	Given $\alpha\in(0,1]$, there exists a constant $~C_{\alpha}$ with the following property.
	For each $\Upsilon>0$ and $q\in\mathbb N$, there exists $\varepsilon_q=\varepsilon_{\Upsilon,q}>0$ such that 
	\begin{itemize}
		\item  for every $\mathcal{C}^{1,\alpha}$ diffeomorphism $f:~M\to M$ satisfying $\|f\|_{\mathcal C^{1,\alpha}}\leq \Upsilon$;
		\item  for every invariant measure $\mu$ of $f$ and every finite partition $\mathcal Q$ with ${\rm Diam}(\mathcal Q)<\varepsilon_q$;
	\end{itemize}
	one has 
	$$h_\mu(f)\le h_{\mu}(f,\cQ)+\frac{1}{\alpha}\big(\frac{1}{q}\int \max_{k}~\log^{+} \|\wedge^k D_xf^q\|~{\rm d}\mu(x)-\lambda^{+}_{\Sigma}(\mu,f)+\frac{1}{q}\big)+\frac{\log(2 q\Upsilon C_{\alpha})}{q}.$$
\end{theoremalph}

	The proof of Theorem \ref{Thm:A-entropy-bound-partition} is based on Burguet's  reparametrization lemma \cite{Bur12}.
	In \cite{Bur12}, Burguet extended the Yomdin-Gromov theory \cite{Gro87,Yom87} to one-dimensional curves and established a remarkable upper bound \cite[Main Proposition]{Bur12} for the ``measure-theoretic tail entropy'' in terms of entropy structures. 
	As an application, he obtained results on symbolic extensions.	
	Burguet’s estimate is formulated for sequences of invariant measures sufficiently close to a given measure. In contrast, our result concerns a single invariant measure and provides a direct estimate relating its metric entropy and partition entropy. 
	This leads to several applications. 
	For instance, see Theorem \ref{Thm:main-perturbation-USC} and Remark \ref{Rem: Important Remark}, it can be used to study perturbations of diffeomorphisms.
	
	The advantage of Burguet's reparametrization lemma is the number of the reparametrizations has better estimates in some sense.
	As a consequence, it has several powerful applications in the study of $C^r$ diffeomorphisms; see, for example, \cite{Bur24,Bur24P}.
	However, obtaining a higher-dimensional version of this result is currently difficult.
	Therefore, we can only assume that $\mu$ has exactly one positive Lyapunov exponent.
    We expect that such a result holds in arbitrary dimensions.
\begin{Conjecture}
	Assume that $\dim M=d>3$. 
	Given $\alpha\in (0,1]$, there exists a constant $~C_{\alpha}$ with the following property.
	For each $\Upsilon>0$ and $q\in\mathbb N$, there exists $\varepsilon_q=\varepsilon_{\Upsilon,q}>0$ such that 
	\begin{itemize}
		\item  for every $\mathcal{C}^{1,\alpha}$ diffeomorphism $f:~M\to M$ satisfying $\|f\|_{\mathcal C^{1,\alpha}}\leq \Upsilon$;
		\item  for every invariant measure $\mu$ of $f$ and every finite partition $\mathcal Q$ with ${\rm Diam}(\mathcal Q)<\varepsilon$;
	\end{itemize}
	one has 
	$$h_\mu(f)\le h_{\mu}(f,\cQ)+\frac{1}{\alpha}\big(\frac{1}{q}\int \max_{1\leq k\leq d}\log^{+}~ \|\wedge^k D_xf^q\|~{\rm d}\mu(x)-\lambda^{+}_{\Sigma}(\mu,f)\big)+\delta(q),~~\lim_{q\to \infty} \delta(q)= 0.$$
\end{Conjecture}
In arbitrary dimensions, it was shown in \cite{LMZ26} that if a diffeomorphism $f$ admits a dominated splitting $E\oplus_{<}F$, then a corresponding result holds for a class of invariant measures.
The estimate obtained there is stronger: there exists $\varepsilon>0$ such that, for every $f$-invariant measure $\mu$ whose Lyapunov exponents along $F$ are all positive and whose Lyapunov exponents along $E$ are all non-positive, and  every finite partition $\mathcal{Q}$ with $\diam(\mathcal Q)<\varepsilon$, we have $h_\mu(f)=h_\mu(f,\mathcal Q)$.

\subsection{Upper semi-continuity of metric entropy and applications}
From the classical work by Ledrappier-Young \cite{LeY85}, it is known that the metric entropy of an ergodic measure depends on the disintegration of the measure along unstable manifolds.
This establishes a fundamental connection between entropy and Lyapunov exponents.
A recent remarkable result by Buzzi-Crovisier-Sarig \cite{BCS22} states that for $\mathcal C^\infty$ surface diffeomorphisms, the continuity of entropy implies the continuity of Lyapunov exponents.
By Theorem \ref{Thm:B-Invar}, we get a opposite direction of \cite{BCS22}: if we have the continuity of Lyapunov exponents of an invariant measure, then we have the upper semi-continuity for the entropy map.

\begin{theoremalph}\label{Thm:main-perturbation-USC}
Assume that $\dim M\le 3$. 
Given any $\mathcal C^{1,\alpha}$ diffeomorphism $f$ and an invariant measure $\mu$ of $f$, 
for any sequence of  $~\mathcal{C}^{1,\alpha}$ diffeomorphisms $\{f_n\}$ and any sequence of probability measures $\{\mu_n\}$ such that $f_n \to f$ in the $\mathcal C^{1,\alpha}$ topology; $\mu_n$ is an invariant measure of $f_n$ for every $n$ and $\mu_n \to \mu$ is the weak star topology.
Then, we have 
$$\limsup_{n\to\infty}h_{\mu_n}(f_n)-h_{\mu}(f)\leq \limsup\limits_{n\to\infty} \frac{1}{\alpha}\big(\lambda_{\Sigma}^+(\mu,f)-\lambda_{\Sigma}^+(\mu_n,f_n)\big).$$
In particular, if $~\lambda_{\Sigma}^+(\mu_n,f_n) \to \lambda_{\Sigma}^+(\mu,f)$, then $\limsup\limits_{n\to\infty}h_{\mu_n}(f_n)\leq h_{\mu}(f)$.
\end{theoremalph}

\begin{Remark}
{\rm 
The theorem is false for $\mathcal C^{1}$ diffeomorphisms.	
Downarowicz-Newhouse \cite[Section 5]{DN05} provide a counterexample: they construct a $\mathcal C^{1}$ surface diffeomorphism $f$ that has a hyperbolic fixed point $p$ and a sequence of ergodic measures $\{ \mu_n \}$ such that $\mu_n \to \delta_p$ (the Dirac measure at $p$), $\lambda^{+}_{\Sigma}(\mu_n) \to \lambda^{+}_{\Sigma}(p)$, but $\lim\limits_{n \to \infty} h_{\mu_n}(f) = \lambda^{+}_{\Sigma}(p) > 0=h_{\delta_p}(f)$.}
\end{Remark} 

\begin{Remark}[Relationship with Theorem 1 of \cite{BL2022}] \label{Rem: Important Remark}
    {\rm Thanks to Burguet's kind reminder, we realized that the case $f_n \equiv f$ in Theorem \ref{Thm:main-perturbation-USC} was already implicitly contained in the paper of Burguet and Liao \cite{BL2022}, and that the case of $\mathcal{C}^{1,\alpha}$ interval maps was proved in \cite{Bur17}.
	
	We first explain how, in the case $f_n\equiv f$, Theorem \ref{Thm:main-perturbation-USC} is related to \cite[Theorem 1]{BL2022}.
	The discussion below is based on personal communications with Burguet.	
	Let $\mu_n$ be a sequence of $f$-invariant measures converging to $\mu$, and let $\{\cP_k\}$ be a sequence of finite partitions such that $\diam(\cP_k)\to 0$ and $\mu(\partial\cP_k)=0$ for every $k$.
	Then, under the additional assumption that $\mu_n$ is ergodic for every $n\geq 1$, \cite[Theorem 1]{BL2022} gives directly
	$$\lim_{k\to \infty} \limsup_{n \rightarrow\infty}  h_{\mu_n}(f)-h_{\mu_n}(f,\cP_k)=\limsup_{n \rightarrow\infty}  h_{\mu_n}(f)-h_{\mu}(f)\leq \frac{1}{\alpha} \limsup_{n \to\infty} \lambda_{\Sigma}^{+}(\mu,f)-\lambda_{\Sigma}^{+}(\mu_n,f).$$
	Indeed, one may identify the entropy structure $\{h_k\}$ in \cite{BL2022} with the sequence of functions $h_k(\nu):=h_\nu(f,\cP_k)$.
	Moreover, the superenvelope function $E(\nu)$ is given by $E(\nu)=\alpha^{-1}\lambda_{\Sigma}^{+}(\nu,f)$.
	
	Note that \cite[Theorem 1]{BL2022} is stated under the assumption that the measures $\mu_n$ close to $\mu$ are ergodic.
	Thus, as stated, it applies directly only to the ergodic case above. 
	We believe that \cite[Theorem 1]{BL2022} suffices to remove the ergodicity assumption on $\mu_n$, but this requires some discussion.
	Indeed, although an invariant measure admits an ergodic decomposition, the fact that an invariant measure $\nu$ is close to $\mu$ does not imply that most of its ergodic components are close to $\mu$.
		
	Theorem \ref{Thm:main-perturbation-USC} provides a perturbative version of this result and does not require the measures $\mu_n$ to be ergodic.
	Its proof is based on Theorem \ref{Thm:B-Invar} and avoids the use of notions such as superenvelopes, entropy structures, fake center foliations, and tail entropy.
	We also provide some applications to the strong positive recurrence (SPR) property and dimension theory.	
	Finally, we believe that the methods developed here may also be applicable to flows generated by four-dimensional $\mathcal{C}^{1,\alpha}$ vector fields, possibly even in the presence of singularities.
}
\end{Remark}	

\subsubsection{Strongly positive recurrence for diffeomorphisms}
More recently, Buzzi-Crovisier-Sarig \cite{BCS25} introduced an important notion of the strongly positive recurrence (SPR) property for diffeomorphisms.  
They proved that SPR diffeomorphisms exhibit exponential mixing and other important statistical properties.
The continuity of Lyapunov exponents is an important property that plays a central role in the study of SPR properties for diffeomorphisms, as shown in \cite{BCS25}. 
The continuity of Lyapunov exponent can be obtained in some natural settings, for instance, the continuity of metric entropy for $\mathcal C^\infty$ surface diffeomorphisms as in \cite{BCS22}.
Our theorem provides some new progress in the study of SPR properties for diffeomorphisms.

\begin{Coro}\label{Thm:SPR-Surface}
	Assume that $f$ is a $\mathcal C^{1,\alpha}$ surface diffeomorphism with $h_{\rm top}(f)>0$.
	Then, $f$ is SPR if and only if for any sequence of ergodic measures $\{\mu_n\}$ with $\mu_n\rightarrow \mu$ and $h_{\mu_n}(f)\rightarrow h_{\rm top}(f)$, one has $\lambda^{+}(\mu_n,f)\rightarrow \lambda^{+}(\mu,f)$.
\end{Coro}

The result can be viewed as a slight improvement of \cite[Theorem B]{BCS25}, since the assumptions on the Lyapunov exponents of $\mu$ are no longer required.
This is because, once the convergence of the Lyapunov exponents is known, Theorem \ref{Thm:main-perturbation-USC} implies that $\mu$ is a measure of maximal entropy. 
Consequently, $\mu$ is $h_{\mathrm{top}}(f)/2$-hyperbolic.

\begin{Coro}\label{Cor:SPR-3}
	Let $f:M\rightarrow M$ be a $\mathcal C^{1,\alpha}$ diffeomorphism on a three-dimensional manifold with positive topological entropy.
	Assume $f$ has the following properties:
	\begin{itemize}
		\item for every sequence of ergodic measures $\{\mu_n\}_{n=1}^{\infty}$ such that $\mu_n$ converging to a measure $\mu$ with $h_{\mu_n}(f)\rightarrow h_{\rm top}(f)$, the measure $\mu$ is hyperbolic and  $\lim\limits_{n \rightarrow\infty} \lambda_{\Sigma}^{+}(\mu_n,f)=\lambda_{\Sigma}^{+}(\mu,f)$.
	\end{itemize}
	Then, we have
	\begin{enumerate}
		\item[(1)] $f$ admits a measure of maximal entropy;
		\item[(2)] there exists $\chi>0$ such that for every ergodic measure of maximal entropy, all its Lyapunov exponents lie outside the interval $[-\chi,\chi]$. 
	\end{enumerate}
	Moreover, for each sequence of ergodic measures $\{\mu_n\}$ with $\mu_n \rightarrow \mu$ and $h_{\mu_n}(f)\rightarrow h_{\rm top}(f)$, there exists $i:=i(\mu)\in \{1,2\}$ such that $\lambda_{i}(x,f)>\chi>-\chi>\lambda_{i+1}(x,f)$ for $\mu$-almost every $x$. 
\end{Coro}
The detailed proof of Corollary \ref{Cor:SPR-3} will be provided in Section \ref{SEC:PofCoE}.
\begin{Remark}\label{re:3-SPR}
	{\rm 
		Based on Corollary~\ref{Cor:SPR-3}, we can improve a bit of \cite[Theorem 3.1]{BCS25} in the following case:
		Assume that $f$ is a $\mathcal C^{1,\alpha}$ diffeomorphism of $M$ ($\dim M=3$)
		with positive topological entropy, then $f$ is SPR if and only if for any sequence of ergodic measures $\{\mu_n\}$ with $\mu_n\rightarrow \mu$ and $h_{\mu_n}(f)\rightarrow h_{\rm top}(f)$, one has $\mu$ is hyperbolic and $\lambda^{+}_{\Sigma}(\mu_n,f)\rightarrow \lambda^{+}_{\Sigma}(\mu,f)$.}
\end{Remark}

\subsubsection{Continuity of the dimension of measures}

Another interesting result is about the Hausdorff dimension of probability measures on closed surface.
Let  $M$ be a $2$-dimensional compact Riemannian manifold without boundary and let $\mu$ be a probability measure on $M$.
The Hausdorff dimension of $\mu$ is defined by
\begin{equation} \label{eq:Y-dimension}
	\dim_H(\mu):= \inf \left\{\dim_H(Z): Z\subset M~\text{with}~\mu(Z)=1 \right\}.
\end{equation}
We have the upper semi-continuity of the Hausdorff dimension under the following conditions.
\begin{Coro}\label{Coro:Dimen-ergodic}
	Let  $M$ be a closed surface, $\{\mu_n\}$ be a sequence of probability measures and $\{f_n\}$ be a sequence of $\mathcal C^{1,\alpha}$ diffeomorphisms such that 
	\begin{itemize}
		\item $\mu_n$ is an ergodic measure of $f_n$ for every $n>0$;
		\item $f_n \rightarrow f$ in the $\mathcal C^{1,\alpha}$ topology;
		\item $\mu_n$ converges to an $f$-ergodic measure $\mu$ and $\lim_{n\to\infty}\limits \lambda^{+}(\mu_n,f_n)=\lambda^{+}(\mu,f)>0$.
	\end{itemize}
	Then, we have~~$\limsup\limits_{n \rightarrow\infty} \dim_H(\mu_n) \leq  \dim_H(\mu)$.
\end{Coro}
The proof of Corollary \ref{Coro:Dimen-ergodic} follows from Theorem \ref{Thm:main-perturbation-USC} and the formula for the Hausdorff dimension of ergodic measures in \cite{Young82}: for every $\mathcal C^{1,\alpha}$ surface diffeomorphism $g$ and every ergodic hyperbolic measure $\nu$ of $g$, one has
$$\dim_H(\nu)=h_{\nu}(g) \left(\frac{1}{\lambda^{+}(\nu,g)}-\frac{1}{\lambda^{-}(\nu,g)}\right).$$
Therefore, under the assumptions of Corollary \ref{Coro:Dimen-ergodic} we have
\begin{align*}
	\limsup_{n \rightarrow\infty} \dim_H(\mu_n)&= \lim_{n \rightarrow\infty} h_{\nu_n}(f_n) \left(\frac{1}{\lambda^{+}(\mu_n,f_n)}-\frac{1}{\lambda^{-}(\mu_n,f_n)}\right) \\
	&= \lim_{n \rightarrow\infty} h_{\nu_n}(f_n) \cdot  \left(\big(\lim_{n \rightarrow\infty} \lambda^{+}(\mu_n,f_n)\big)^{-1}-\big(\lim_{n \rightarrow\infty} \lambda^{-}(\mu_n,f_n)\big)^{-1}\right) \\
	&\leq h_{\mu}(f) \left(\frac{1}{\lambda^{+}(\mu,f)}-\frac{1}{\lambda^{-}(\mu,f)}\right)=\dim_H(\mu).
\end{align*}
This completes the proof of Corollary \ref{Coro:Dimen-ergodic}.
\begin{Coro}\label{Cor:Dim-MME}
	Let $f$ be a $\mathcal{C}^{\infty}$ surface diffeomorphism with $h_{\rm top}(f)>0$. 
	Then, for every sequence of ergodic measures $\{\mu_n\}$ with $\mu_n \rightarrow \mu$ and 
	$h_{\mu_n}(f) \rightarrow h_{\rm top}(f)$ as $n\rightarrow \infty$, we have that $\lim\limits_{n\rightarrow \infty} \dim_H(\mu_n)=\dim_H(\mu)$.
\end{Coro}

We now explain how this result is obtained.
In the setting of Corollary \ref{Cor:Dim-MME}, by \cite[Theorem B]{BLY25} the limit measure $\mu$ is an ergodic measure of maximal entropy. 
Moreover, by \cite{Bur24P} we have $\lambda^{+}(\mu_n,f) \rightarrow \lambda^{+}(\mu,f)$ as $n\rightarrow \infty$.  
Therefore, Corollary \ref{Cor:Dim-MME} follows.

\begin{Remark}
	By the main Theorem of \cite{Bur24P} and \cite[Theorem B]{BLY25}.
	This result remains valid in the $\mathcal{C}^{r,\alpha}$ setting, provided that $h_{\rm top}(f)>(r+\alpha)^{-1} \lim\limits_{n \to \infty} \tfrac{1}{n} \log  \min\{\|Df^{n}\|_{0},~\|Df^{-n}\|_{0}\}$.
\end{Remark}

\section{Entropy estimate for ergodic measures with finite partitions}

In this section, we give the proofs of Theorems \ref{Thm:A-entropy-bound-partition}.
For the main theorems in this paper, we will only consider low regularity, i.e., the $\mathcal C^{1,\alpha}$ case. 
However, Theorem \ref{Thm:A-entropy-bound-partition} may have general interest when the map is $\mathcal C^{r,\alpha}$ for some large $r>1$.
We begin by stating a result under stronger differentiability assumptions (the notations will be introduced in the next section).
\begin{Theorem}\label{Thm:entropy-bound-partition}
Given $r\in\mathbb N$ and $\alpha\in[0,1]$ satisfying $r+\alpha>1$, there exists a constant $C_{r,\alpha}$ with the following property.
For each $\Upsilon>0$ and $q\in\mathbb N$, there exists $\varepsilon=\varepsilon_{\Upsilon,q}>0$ such that 
\begin{itemize}
	\item  for every $\mathcal C^{r,\alpha}$ diffeomorphism $f:~M\to M$ satisfying $\|f\|_{\mathcal C^{r,\alpha}}\leq \Upsilon$;
	\item  for every ergodic measure $\mu$ of $f$ with exactly one positive Lyapunov exponent;
	\item  for every finite partition $\mathcal Q$ with ${\rm Diam}(\mathcal Q)<\varepsilon$;
\end{itemize}
one has 
$$h_\mu(f)\le h_\mu(f,\mathcal Q)+\frac{1}{r-1+\alpha}\big[\frac{1}{q}\int\log\|D_xf^q\|{\rm d}\mu-\lambda^+(\mu,f)+\frac{1}{q}\big]+\frac{\log(2q \Upsilon\cdot C_{r,\alpha})}{q}.$$
\end{Theorem}

\begin{Remark}\label{Rem:Zerobou}
	{\rm Note that Theorem \ref{Thm:A-entropy-bound-partition} follows from Theorem \ref{Thm:entropy-bound-partition}. 
	Without loss of generality, we can assume that $\mu(\partial Q)=0$.
	Indeed, for any finite partition $\mathcal{Q}$ with ${\rm Diam}(\mathcal Q)<\varepsilon$ and every $\eta>0$, one can always find a finite partition $\cQ'$ satisfying $\diam(\cQ')<\varepsilon$ and $\mu(\partial \cQ')=0$ such that $|h_{\mu}(f,\cQ)-h_{\mu}(f,\cQ')|< \eta$.
	Consequently, if Theorem \ref{Thm:entropy-bound-partition} holds for every partition $\mathcal Q'$ with $\diam(\mathcal Q')<\varepsilon$ and $\mu(\partial\mathcal Q')=0$, then it also holds for every partition $\mathcal Q$ with $\diam(\mathcal Q)<\varepsilon$.}
\end{Remark}

Recall that $\lambda^-(\mu,f)=-\lambda^+(\mu,f^{-1})$.
By considering the $\mathcal C^{r,\alpha}$ diffeomorphism $f^{-1}$, one obtains the following symmetric version of Theorem~\ref{Thm:entropy-bound-partition}. 

\begin{Theorem}\label{Thm:entropy-bound-partition-inverse}
Given $r\in\mathbb N$ and $\alpha\in[0,1]$ satisfying $r+\alpha>1$, there exists a constant $C_{r,\alpha}$ with the following property.
For each $\Upsilon>0$ and $q\in\mathbb N$, there exists $\varepsilon=\varepsilon_{\Upsilon,q}>0$ such that 
\begin{itemize}
	\item  for every $\mathcal C^{r,\alpha}$ diffeomorphism $f:~M\to M$ satisfying $\|f\|_{\mathcal C^{r,\alpha}}\leq \Upsilon$;
	\item  for every ergodic measure $\mu$ of $f$ with exactly one negative Lyapunov exponent;
	\item  for every finite partition $\mathcal Q$ with ${\rm Diam}(\mathcal Q)<\varepsilon$;
\end{itemize}
one has 
$$h_\mu(f)\le h_\mu(f,\mathcal Q)+\frac{1}{r-1+\alpha}\big[\frac{1}{q}\int\log\|D_xf^{-q}\|{\rm d}\mu+\lambda^{-}(\mu,f)+\frac{1}{q}\big]+\frac{\log(2q \Upsilon\cdot C_{r,\alpha})}{q}.$$
\end{Theorem}

\begin{Remark}
	{\rm
		Tail entropy, or local entropy, which introduced by Buzzi \cite{Buz97} and  Newhouse \cite{New89}, is used to estimate the degree to which entropy fails to be upper semi-continuous.   
		For every three-dimensional diffeomorphism, any ergodic measure  has either exactly one positive Lyapunov exponent or exactly one negative Lyapunov exponent. 
		Then, in a sense, Theorem \ref{Thm:entropy-bound-partition} and Theorem \ref{Thm:entropy-bound-partition-inverse} provide new bounds on local entropy with respect to measures. }
\end{Remark}

In the remainder of this section, we will provide a detailed proof of Theorem \ref{Thm:entropy-bound-partition}.
The proof of Theorem \ref{Thm:entropy-bound-partition-inverse} can be proved similarly by considering the $\mathcal C^{r,\alpha}$ diffeomorphism $f^{-1}$.
By Remark \ref{Rem:Zerobou}, we can assume that $\mu(\partial \cQ)=0$.

\subsection{Preliminaries}
Assume that $M$ is a compact boundaryless Riemannian manifold.
For a $\mathcal C^{r,\alpha}$ diffeomorphism $f:M\to M$ with $r\in\mathbb N$ and $\alpha\in[0,1]$, we mean
\begin{itemize}
	\item if $\alpha=0$, it is just the usual $\mathcal C^r$ diffeomorphism with $r\in\mathbb N$;
	\item if $\alpha=1$, it is a $\mathcal C^r$ diffeomorphism, and its $\mathcal C^r$ derivative $D^r f$ is a Lipschitz map;
	\item if $\alpha\in(0,1)$,  it is a $\mathcal  C^r$ diffeomorphism, and its $\mathcal C^r$ derivative $D^r f$ is a $\alpha$-H\"older map.
\end{itemize} 

Assume that $X$ is a compact metric space and $Y$ is a Banach space.
For $\alpha\in (0,1]$ and a $\alpha$-H\"{o}lder continuous map $H:X\to Y$, define
\begin{equation}\label{eq:normd}
	\|H\|_0=\sup_{x\in X}\|H(x)\|,~~\|H\|_\alpha=\sup \left \{ \dfrac{d(H(x),H(y))}{d(x,y)^\alpha}:x\neq y,~x,y\in X \right\}.
\end{equation}
For $\mathcal C^{r,\alpha}$ diffeomorphism $f:~M\to M$, define
$$\|f\|_{\mathcal C^{r,\alpha}}:=\max\Big\{ \max_{1\le j\le r}\|D^jf\|_0,~\max_{1\le j\le r}\|D^jf^{-1}\|_0,~\|D^rf\|_\alpha,~\|D^rf^{-1}\|_\alpha \Big\}.$$

We recall some fundamental definition and properties of entropies.
Let $\mu$ be a  probability measure.
Given a finite partition $\mathcal P$, define the static entropy function
$$H_\mu(\mathcal P):=\sum_{P\in\mathcal P}-\mu(P)\log\mu(P).$$
By definition, one has 
\begin{equation}\label{eq:estimate-conditional}
	H_\mu(\mathcal P)\le \log\#\{P\in\mathcal P:~\mu(P)>0\}.
\end{equation}
Given two finite partitions $\mathcal P$ and $\mathcal Q=\{Q_1,\cdots,Q_k\}$, define the conditional entropy
$$H_\mu(\mathcal P\big|\mathcal Q):=\sum_{j=1}^k\mu(Q_j)H_{\mu_j}(\mathcal P),$$
where $\mu_j(\cdot)=\frac{\mu(Q_j\cap \cdot)}{\mu(Q_j)}$ denotes the normalization of $\mu$ restricted on $Q_j$.

For an $f$-invariant measure $\mu$ and a finite partition $\mathcal P$, denote by 
$$\mathcal P^n:=\mathcal P^{n,f}=\bigvee_{j=0}^{n-1}f^{-j}(\mathcal P).$$
The metric entropy of $\mu$ associated to a partition $\mathcal P$ is defined as 
$$h_\mu(f,\mathcal P):=\lim_{n\to\infty}\frac{1}{n}H_\mu(\mathcal P^n);$$
and the metric entropy of $\mu$ is defined as
$$h_\mu(f):=\sup \{h_\mu(f,\mathcal P): \mathcal P~\textrm{is a finite partition}\}.$$
Note that $h_\mu(f)=h_\mu(f^{-1})$ and $h_\mu(f,\mathcal P)=h_\mu(f^{-1},\mathcal P)$ for every finite partition $\mathcal{P}$ (see \cite[Theorem 4.13]{Wal82}).

\subsection{The Reparametrization Lemma and the choice of $C_{r,\alpha}$}
We recall Burguet's reparametrization lemma \cite{Bur24} for $\mathcal C^{r,\alpha}$ diffeomorphisms. 
A $\mathcal C^{r,\alpha}$ curve $\sigma:[-1,1]\rightarrow M$ with $r+\alpha>1$ is said to be \emph{bounded}, if it satisfies the following conditions
\begin{itemize}
\item if $r\ge 2$, then 
$$\sup_{2\leq s\leq r} \|D^s \sigma\|_{0} \leq \frac{1}{6} \|D\sigma\|_{0},~~~ \|D^r \sigma\|_{\alpha} \leq \frac{1}{6} \|D\sigma\|_{0}.$$
\item if $r=1$ and $\alpha\in(0,1]$, then 
$$ \|D \sigma\|_{\alpha} \leq \frac{1}{6} \|D\sigma\|_{0}.$$
\end{itemize}
A bounded curve $\sigma:[-1,1]\rightarrow M$ is said to be \emph{strongly $\varepsilon$-bounded}, if $\|D\sigma\|_{0}\leq \varepsilon$.
For a curve  $\sigma:[-1,1]\rightarrow M$, denote by $\sigma_*=\sigma([-1,1])$ the image of $\sigma$.

\begin{Lemma}[\cite{Bur24P}, Lemma 12]\label{Lem:local-reparametrization}
Given $r\in\mathbb N$ and $\alpha\in[0,1]$ satisfying $r+\alpha>1$, there exists a constant $C_{r,\alpha}$ with the following property.
Given $\Omega>0$, there exists $\varepsilon_\Omega>0$ such that if $g$ is a $\mathcal C^{r,\alpha}$ diffeomorphism with 
\begin{equation}\label{e.Omega-bounded}
\max_{1\le j\le r}\|D^jg\|_0<\Omega,~~~\|D^rg\|_\alpha<\Omega
\end{equation}
 then for any strongly $\varepsilon$-bounded $\mathcal C^{r,\alpha}$ curve $\sigma:~[-1,1]\to M$ with $\varepsilon\in(0,\varepsilon_\Omega)$ and any $\chi^+,\chi \in\mathbb Z$, there is a family of affine reparametrizations $\Theta$ such that
\begin{enumerate}
\item[(1)] $\{t\in[-1,1]:\lceil \log \| D_{\sigma(t)}g \| \rceil=\chi^+,~\lceil\log\|D_{\sigma(t)}g|_{T_{\sigma(t)} \sigma_*}\|\rceil=\chi\}\subset \bigcup_{\theta\in\Theta}\theta([-1,1])$;
\item[(2)] $g\circ \sigma\circ\theta$ is bounded for any $\theta\in\Theta$;
\item[(3)] $\#\Theta\le C_{r,\alpha} \exp(\frac{\chi^+-\chi}{r+\alpha-1})$;
\end{enumerate}
where $\lceil a \rceil$ denotes the smallest integer that is larger than or equal to $a$.
\end{Lemma}

This $\mathcal C^{r,\alpha}$ version of the reparametrization lemma is parallel to Burguet's work, which consider the case of $\alpha=0$ in \cite{Bur24P}. However, some preparations for the case $\alpha\in(0,1]$ were previously carried out in \cite{Bur12}. 
For completeness, we provide a detailed proof in the appendix.

The constant $C_{r,\alpha}$ appearing in the statement of Theorem~\ref{Thm:entropy-bound-partition} is precisely the one chosen from Lemma~\ref{Lem:local-reparametrization}.

\subsection{Choice of $\varepsilon:=\varepsilon_{\Upsilon,q}$} \label{SEC:Vepsilon}
Since $M$ is compact, we can choose $r(M)>0$ such that $\exp_x^{-1}:B(x,2r(M))\rightarrow T_x M$ is a $\mathcal C^{\infty}$ embedding. 
Then, by changing the metric if necessary, for each bounded curve $\sigma:[-1,1] \rightarrow  M$ with ${\rm diam}(\sigma_{\ast})<r(M)$, for any $y\in M$ and any $\varepsilon>0$, if $\sigma_{\ast}\cap B(y,\varepsilon)\neq \emptyset$, then we can choose a reparametrization $\theta$ such that $(\sigma\circ \theta)_{\ast}=\sigma_{\ast}\cap B(y,\varepsilon)$ and $\sigma\circ \theta$ is strongly $2\varepsilon$-bounded.

Given $\Upsilon>0$ and $q\in\mathbb N$, there exists $\Omega>0$ with the following properties:  
if $f:~M\to M$ is a $\mathcal C^{r,\alpha}$ diffeomorphism  satisfying
$$\max_{1\le j\le r}\|D^jf\|_0<\Upsilon,~\max_{1\le j\le r}\|D^jf^{-1}\|_0<\Upsilon,~\|D^rf\|_\alpha<\Upsilon,~\|D^rf^{-1}\|_\alpha<\Upsilon,$$
then for any $1\le k\le q$, one has 
$$\max_{1\le j\le r}\|D^jf^k\|_0<\Omega,~\max_{1\le j\le r}\|D^jf^{-k}\|_0<\Omega,~\|D^rf^k\|_\alpha<\Omega,~\|D^rf^{-k}\|_\alpha<\Omega.$$
Now we choose $\varepsilon_\Omega$ from Lemma~\ref{Lem:local-reparametrization}. We then choose $\varepsilon:=\varepsilon_{\Upsilon,q}\in (0, \frac{\varepsilon_\Omega}{4})$ such that
\begin{equation}\label{e.choose-varepsilon}
	2(\Omega+2)\varepsilon<\min\{1,r(M)\}. 
\end{equation}

\subsection{Choose a set $K$, a strongly $\varepsilon$-bounded curve $\sigma$ and a partition $\mathcal P$}
Fix an ergodic measure $\mu$ and fix the partition $\mathcal Q$ satisfying the conditions in Theorem~\ref{Thm:entropy-bound-partition}, i.e., 
$\mu$ has exactly one positive Lyapunov exponent, ${\rm Diam}(\mathcal Q)<\varepsilon$ and $\mu(\partial\mathcal Q)=0$.

Since $\mu$ has exactly one positive Lyapunov exponent, it follows from \cite{LeS82} that there exists a measurable partition $\xi$ subordinate to the one-dimensional Pesin unstable foliation $W^u$. 
By the Rokhlin disintegration theorem \cite{Rok67}, we denote by $\{ \mu_{\xi(x)} \}$ the family of conditional measures of $\mu$ with respect to the measurable partition $\xi$. 

The following proposition is a consequence of Ledrappier-Young's result \cite{LeY85}, for the proof, see \cite[Proposition 2.1, Proposition 2.2]{LuY24}.
\begin{Proposition}\label{Prop:Two Balls}
	For every $\tau>0$ and every $\delta\in(0,1)$, there exists $K\subset M$ with $\mu(K)>1-\delta$ and $\rho:=\rho_K>0$, such that for every $x\in K$, every measurable set $\Sigma \subset W^{u}_{{\rm loc}}(x)$ with $\mu_{\xi(x)}(\Sigma\cap K)>0$, and every finite partition $\cP$ with $ {\rm Diam}(\cP)<\rho$, one has
	\begin{equation}\label{eq:LimK}
		h_{\mu}(f)\le\liminf_{n\to\infty}\frac{1}{n} H_{\mu_{\xi(x),K,\Sigma}} (\cP^n)+\tau,
	\end{equation}
	where $\mu_{\xi(x),K,\Sigma}(\cdot):=\frac{\mu_{\xi(x)}(K\cap \Sigma \cap \cdot )}{\mu_{\xi(x)} (K\cap \Sigma)}$.
\end{Proposition}
Given an auxiliary constant $\tau>0$. We choose a compact set $K$ with the following properties:
\begin{itemize}
\item $\mu(K)>\frac{1}{2}$ and $K$ satisfies the conclusion of Proposition \ref{Prop:Two Balls};  
\item the following convergences hold uniformly for $x\in K$
\begin{equation}\label{e.uniform-converge}
\frac{1}{n}\sum_{j=0}^{n-1}\delta_{f^{j}(x)}\to \mu,~~\frac{1}{n} \log \|D_xf^{n}|_{E^{u}(x)}\| \rightarrow \lambda^{+}(\mu,f),
\end{equation}
where $\delta_x$ denotes the Dirac measure at $x$ and $E^u$ is the one-dimensional measurable bundle associated to the positive Lyapunov exponent;
\item for every $c\in \{0,\cdots,q-1\}$, the following convergence holds uniformly for $x\in K$
$$ \lim_{m \to \infty} \frac{1}{m} \sum_{j=0}^{m-1} \log \|D_{f^{qj+c}(x)}f^q\|:=\phi_{c}(x),$$
where $\phi_{c}:M\rightarrow \RR$ is an $f^q$-invariant  measurable function and for every $x\in K$ one has
$\sum_{c=0}^{q-1} \phi_{c}(x)=\int \log \|D_{x}f^q\| {\rm d}\mu$.

\end{itemize}

Choose a point $x_0\in K$ and a strongly $\varepsilon$-bounded curve $\sigma:[-1,1]\rightarrow W^{u}_{\loc}(x_0)$ such that $\mu_{\xi(x_0)}(\sigma_{\ast}\cap K)>0$. 
Consider the measure $\mu_{\sigma}$ defined as follows
$$\mu_{\sigma}(A):=\frac{\mu_{\xi(x_0)}(K\cap \sigma_{\ast} \cap A)}{\mu_{\xi(x_0)}(K\cap \sigma_{\ast} )},~~\forall~\textrm{Borel set}~A.$$
By Proposition \ref{Prop:Two Balls}, there exists a finite partition $\mathcal{P}$ with $\mu(\partial \mathcal P)=0$ such that
$$h_{\mu}(f) \le \liminf_{n\to\infty}\frac{1}{n} H_{\mu_{\sigma}}(\mathcal P^n)+\tau.$$
Using the properties of conditional entropy, one has that
$$H_{\mu_{\sigma}}(\mathcal P^n)\le  H_{\mu_{\sigma}}(\mathcal Q^n)+H_{\mu_{\sigma}}(\mathcal P^n\big| \mathcal Q^n),$$
Therefore, we obtain
\begin{equation}\label{eq: Entropymu}
h_{\mu}(f) \le \limsup_{n\to\infty}\frac{1}{n}H_{\mu_{\sigma}}(\mathcal Q^n)+\limsup_{n\to\infty}\frac{1}{n}H_{\mu_{\sigma}}(\mathcal P^n\big| \mathcal Q^n)+\tau.
\end{equation}

\subsection{Estimate for $H_{\mu_{\sigma}}(\mathcal Q^n)$}
Given an integer $N \in\mathbb N$, by the concave property of the function $H$ (see \cite[Section 8.2]{Wal82}), for every $n>N$ one has

$$\frac{1}{n}H_{\mu_{\sigma}}(\mathcal Q^n)\le \frac{1}{N}H_{\frac{1}{n}\sum_{i=0}^{n-1}f_*^i \mu_{\sigma}}(\mathcal Q^N)+\frac{2N}{n}\log\#\mathcal Q.$$
By the definition of $\mu_{\sigma}$, one has
$$\mu_{\sigma}=\frac{1}{\mu_{\xi(x_0)}(K\cap \sigma_*)}\int_{K\cap \sigma_*} \delta_x~{\rm d}\mu_{\xi(x_0)}(x),$$
By the choice of $K$, (see \eqref{e.uniform-converge}), one has that
$$\lim_{n\to\infty}\frac{1}{n}\sum_{i=0}^{n-1}f_*^i \mu_{\sigma}=\mu.$$
Since $\mu(\partial\mathcal Q^N)=0$ for every $N>0$ (which can be deduced from $\mu(\partial\mathcal{Q})=0$), one has 
$$\lim_{n\to\infty}H_{\frac{1}{n}\sum_{i=0}^{n-1}f_*^i \mu_{\sigma}}(\mathcal Q^N)=H_\mu(\mathcal Q^N).$$
Thus, for every $N\in \NN$ we have 
\begin{equation}
	\limsup_{n\to\infty}\frac{1}{n} H_{\mu_{\sigma}}(\mathcal Q^n)\le \frac{1}{N}H_\mu(\mathcal Q^N).
\end{equation}
Consequently, we have 
\begin{equation}\label{eq:firstpart}
	\limsup_{n\to\infty}\frac{1}{n} H_{\mu_{\sigma}}(\mathcal Q^n)\le h_{\mu}(f,\mathcal{Q}).
\end{equation}

\subsection{The estimate of $H_{\mu_{\sigma}}(\mathcal P^n\big| \mathcal Q^n)$}
Recall that $\sigma:~[-1,1]\to W^u_{\rm loc}(x_{0})$ is a strongly $\varepsilon$-bounded curve, where $\varepsilon$ is chosen in Section \ref{SEC:Vepsilon}.

\begin{Proposition}\label{Pro:cover-n-step}
	For every $y\in \sigma_*\cap K$ and every $n\in \NN$, there exists a family of reparametrizations $\Gamma_n$ such that
\begin{enumerate}
\item[(1)] $\sigma_*\cap K\cap B_n(y,\varepsilon)\subset \bigcup_{\gamma\in\Gamma_n}\sigma\circ\gamma([-1,1])$;
\item[(2)] for any $0\le j\le n-1$, $\|Df^j\circ\sigma\circ\gamma\|_{0}\le 1$;
\item[(3)] $\limsup\limits_{n\to\infty}\frac{1}{n}\log\#\Gamma_n\le \frac{1}{r-1+\alpha}\big(\frac{1}{q}\int \log\|Df^q\|{\rm d}\mu-\lambda^+(\mu,f)+\frac{1}{q}\big)+\frac{\log(2q \Upsilon \cdot C_{r,\alpha})}{q}$.
\end{enumerate}

\end{Proposition}
\begin{proof}
	By the choice of $K$, we have
	$$\sum_{k=0}^{q-1} \lim_{n\rightarrow \infty}\frac{1}{n} \sum_{j=0}^{n-1} \log \|D_{f^{jq+k}(x)} f^q\|=\lim_{n\rightarrow \infty}\frac{1}{n} \sum_{j=0}^{n-1} \log \|D_{f^j(x)} f^q\| =\int \log \|D_xf^q \| {\rm d} \mu(x).$$
	Hence, for every $x\in K$ there exists $c(x)\in \{0,\cdots,q-1\}$ such that 
	\begin{equation}\label{e.choose-constant-c}
		\lim_{m \rightarrow\infty}\frac{1}{m} \sum_{j=0}^{m-1} \log \|D_{f^{jq+c(x)}(x)} f^q\|\leq \frac{1}{q}\int \log \|D_xf^q \| {\rm d} \mu(x),
	\end{equation}
	where $\lfloor a \rfloor$ is the largest integer less than or equal to $a$.
We decompose $K$ to be the union of $\{K_c\}_{c=0}^{q-1}$ such that for any $x\in K_c$, one has that $c(x)=c$. 

We now fix one $K_c$ such that $\sigma_*\cap K_c\cap B_n(y,\varepsilon)\neq \emptyset$, and we take $m=\lfloor (n-c)/q \rfloor$.
We decompose $\sigma_*\cap K_c\cap B_n(y,\varepsilon)$ into subsets $\Sigma((\chi_j^+,\chi_j)_{j=0}^{m})$, where the set $\Sigma((\chi_j^+,\chi_j)_{j=0}^{m	})$ is defined as the points $z\in \sigma_*\cap K_c\cap B_n(y,\varepsilon)$ such that
\begin{itemize}
\item for $j=0$,
$$\lceil\log\|D_zf^c\|\rceil=\chi^+_0,~\lceil\log\|D_zf^{c}|_{T_{z}(\sigma_*)}\|\rceil=\chi_0;$$
\item for any $1\le j\le m$,
$$\lceil\log\|D_{f^{q(j-1)+c}(z)}f^q\|\rceil=\chi^+_j,~\lceil\log\|D_{f^{q(j-1)+c}(z)}f^{q}|_{T_{f^{q(j-1)+c}(z)}{(f^{q(j-1)+c}(\sigma_*))}}\|\rceil=\chi_j$$
\end{itemize}

Note that a bounded curve constrained within a small ball of radius  $\varepsilon$  is strongly  $2\varepsilon$-bounded.
Indeed, let $\gamma:[-1,1]\to M$ be a bounded curve. 
For every $x_{\ast}\in M$ such that $\gamma_{\ast}\cap B(x_{\ast},\varepsilon)$, consider the interval
$$[a,b]:=\big[\min\{t\in[-1,1]:\gamma(t)\in B(x_{\ast},\varepsilon)\},~\max\{t\in[-1,1]:\gamma(t)\in B(x_{\ast},\varepsilon)\}\big]$$
and let $\theta:[-1,1]\to [a,b]$ be the affine reparametrization.
Then, $\gamma \circ \theta$ is strongly $2\varepsilon$-bounded, since the length of $\gamma \circ \theta$ is at most $4\varepsilon$.  

For each set $\Sigma((\chi_j^+,\chi_j)_{j=0}^{m})$, one applies Lemma~\ref{Lem:local-reparametrization} first for $f^c$ and then for $f^q$ inductively $m$ times. 
For $y\in \sigma_{\ast} \cap K$, the set 
\begin{equation}\label{e.non-emptyset}
\{z\in \Sigma((\chi_j^+,\chi_j)_{j=0}^{m-1}):~d(f^{c+jq}(z),f^{c+jq}(y))<\varepsilon,~\forall 0\le j< m\}
\end{equation}
 admits a family of reparametrizations $\Gamma((\chi_{j}^{+},\chi_j)_{j=0}^{m})$ such that
\begin{itemize}
\item $f^{jq+c}\circ \sigma\circ\theta$ is strongly $2\varepsilon$-bounded for any $\theta\in\Gamma((\chi_j^+,\chi_j)_{j=0}^{m})$ and any $0\le j<m$;
\item $\#\Gamma((\chi_j^+,\chi_j)_{j=0}^{m})\le C_{r,\alpha}^{m+1} \exp(\frac{1}{r-1+\alpha}\sum_{j=0}^{m}(\chi_j^+-\chi_j))$.
\end{itemize}
Since $\max\{\chi_j^{+},\chi_j\}\leq \log \Upsilon$ for every $0 \leq j \leq m$, there are at most $(q \log \Upsilon)^{2m+2}$ possible choice of $(\chi_j^+,\chi_j)_{j=0}^{m}$ such that the set in \eqref{e.non-emptyset} is non-empty. See \cite{Bur24} for instance. 

By the definitions, for every $z\in \Sigma((\chi_j^+,\chi_j)_{j=0}^{m-1})\cap K_{c}$ one has 
$$
\sum_{j = 0}^{m} \chi_{j}^{+} \leq \log \| D_zf^{c}\|+\sum_{j = 0}^{m-1}\log\| D_{f^{jq+c}(z)}f^{q}\|+m + 1
$$
Thus, by \eqref{e.choose-constant-c} one has 
$$\limsup_{n\rightarrow\infty}\frac{1}{n}\sum_{j = 0}^{m}\chi_{j}^{+} \leq \lim_{n\rightarrow\infty}\frac{1}{n}\left(\sum_{j =0}^{m-1}\log\| D_{f^{jq+ c}(z)}f^{q}\|+m\right)\leq\frac{1}{q}\int\log\| D_xf^{q}\| {\rm d}\mu(x)+\frac{1}{q}.$$
Since
$$\sum_{j = 0}^{m} \chi_{j} \geq \log\| D_zf^{i}\big|_{E^{u}(z)}\|+\sum_{j = 0}^{m - 1}\log\| Df^{q}\big|_{E^{u}(f^{jq+ i}(z))}\|=\log \|D_zf^{mq+c}|_{E^u(z)}\|,$$
one has that
$$\liminf_{n\rightarrow\infty}\frac{1}{n}\sum_{j = 0}^{m} \chi_{j} \geq \lim_{n\rightarrow\infty}\frac{1}{n}\log\| Df^{n}\big|_{E^{u}(z)}\|=\lambda^{+}(\mu,f).$$

Take $\Gamma_n$ be the union of all these possible $\Gamma((\chi_j^+,\chi_j)_{j=0}^{m})$ for all $K_c$. 
By the construction, Item (1) holds.

For any $\gamma\in \Gamma_n$, there exists $0\leq c <q $ such that $f^{jq+c}\circ \sigma\circ\theta$ is strongly $2\varepsilon$-bounded for any $0\le j \le \lfloor (n-c)/q\rfloor$.
By the choice of $\varepsilon$ (Equation~\eqref{e.choose-varepsilon}), one has that $\|Df^i \circ\sigma\circ\gamma\|\le 1$ for every $0\leq i<n$. 
Thus, Item~(2) holds. 
It remains to estimate the cardinality of $\Gamma_n$. 

By the construction of $\Gamma_n$, one has that
\begin{align*}
	&~\limsup_{n\to\infty}\frac{1}{n}\log\#\Gamma_n\\
	\le&~\limsup_{n\to\infty}\frac{1}{n}\log q \cdot \left(( q \log \Upsilon )^2C_{r,\alpha} \right)^{m+1}+\frac{1}{r+\alpha-1}\left(\frac{1}{q}\int \log \|Df^q\| {\rm d}\mu- \lambda^{+}(\mu,f)+\frac{1}{q}\right)
\end{align*}
Thus, Item (3) can be concluded:
$$\limsup_{n\to\infty}\frac{1}{n}\log\#\Gamma_n\le \frac{1}{r-1+\alpha}\left(\frac{1}{q}\int \log\|Df^q\|{\rm d}\mu-\lambda^+(\mu,f)+\frac{1}{q}\right)+\frac{\log(2q\Upsilon C_{r,\alpha} )}{q}.$$
This completes the proof of Proposition \ref{Pro:cover-n-step}.
\end{proof}

Since the diameter of $\mathcal Q$ is smaller than $\varepsilon$, from Proposition~\ref{Pro:cover-n-step}, one has the following corollary directly. 

\begin{Corollary}\label{Cor:cover-n-step-partition}
For every $n\in \NN$ and every $Q_n\in \mathcal Q^n$, there is a family of reparametrizations $\Gamma(Q_n)$, such that
\begin{enumerate}
\item[(1)] $\sigma_*\cap K\cap Q_n\subset \bigcup_{\gamma\in\Gamma_n}\sigma\circ\gamma([-1,1])$;
\item[(2)] $\|Df^j\circ\sigma\circ\gamma\|\le 1$ for any $0\le j\le n-1$;
\item[(3)] $\limsup\limits_{n\to\infty}\dfrac{1}{n}\sup\limits_{Q_n\in \mathcal{Q}^n}\log\#\Gamma(Q_n)\le \frac{1}{r-1+\alpha}\big(\frac{1}{q}\int \log\|Df^q\|{\rm d}\mu-\lambda^{+}(\mu,f)+\frac{1}{q}\big)+\frac{\log(2q\Upsilon C_{r,\alpha})}{q}$.
\end{enumerate}

\end{Corollary}

\begin{Theorem}\label{Thm:conditional-entropy}
$$\limsup_{n\to\infty}\frac{1}{n}H_{\mu_{\sigma}}(\mathcal P^n\big| \mathcal Q^n)\le \frac{1}{r-1+\alpha}\big(\frac{1}{q}\int \log\|Df^q\|{\rm d}\mu-\lambda^+(\mu,f)+\frac{1}{q}\big)+\frac{\log(2q\Upsilon C_{r,\alpha})}{q}.$$

\end{Theorem}
\begin{proof}
By the definition of the conditional entropy
$$H_{\mu_{\sigma}}(\mathcal P^n\big| \mathcal Q^n)=\sum_{Q\in \mathcal Q^n:~\mu_{\sigma}(Q)>0} \mu_{\sigma}(Q_n) \cdot H_{\mu_{\sigma,Q_n}}(\mathcal P^n),$$
where $\mu_{\sigma,Q_n}=\frac{\mu_{\sigma}(Q_n \cap \cdot )}{\mu_{\sigma}(Q_n)}$ is the normalization of $\mu_{\sigma}$ on $Q_n$.

By \eqref{eq:estimate-conditional}, one has that 
\begin{align*}
H_{\mu_{\sigma}}(\mathcal P^n | \mathcal Q^n)&\le \sum_{Q_n\in \mathcal Q^n,~\mu_{\sigma}(Q_n)>0} \mu_{\sigma}(Q_n) \cdot \log \#\{P_n\in \mathcal P^n:~P_n\cap \sigma_*\cap K\cap Q_n\neq\emptyset\}.
\end{align*}
By Corollary~\ref{Cor:cover-n-step-partition}, for each $n>0$ and each $Q_n \in \mathcal{Q}^n$, there is a reparametrization family $\Gamma(Q_n)$ such that 
\begin{enumerate}
\item $\sigma_*\cap K\cap Q_n\subset \bigcup_{\gamma\in\Gamma_n}\sigma\circ\gamma([-1,1])$;
\item $\|Df^j\circ\sigma\circ\gamma\|_{0} \le 1$  for every $0\le j\le n-1$;
\item $\limsup\limits_{n\to\infty}\dfrac{1}{n}\sup\limits_{Q_n\in \mathcal{Q}^n}\log\#\Gamma(Q_n)\le  \frac{1}{r-1+\alpha}\big(\frac{1}{q}\int \log\|Df^q\|{\rm d}\mu-\lambda^{+}(\mu,f)+\frac{1}{q}\big)+\frac{\log(2q\Upsilon C_{r,\alpha})}{q}$.
\end{enumerate}
Since we have
$$\#\{P_n\in \mathcal P^n:~P_n\cap \sigma_*\cap K\cap Q_n\neq\emptyset\}\le \sum_{\gamma\in\Gamma(Q_n)}\#\{P_n\in \mathcal P^n:~P_n\cap K\cap ( \sigma\circ\gamma)_*\neq\emptyset\},$$
it suffices to estimate $\#\{P_n\in \mathcal P^n:~P_n\cap K \cap(\sigma\circ\gamma)_*\neq\emptyset\}$ for each $\gamma\in \Gamma_n$.
Recall that we assume $\mu(\partial \mathcal{P})=0$, then inspired by \cite[Page 1498]{Bur24P}, see also \cite[Proposition 5.2]{LuY24}, one has that
$$\limsup_{n\to\infty}\frac{1}{n}\log\sup_{\gamma\in\Gamma(Q_n)}\#\{P\in \mathcal P^n:~P\cap (\sigma\circ\gamma)_*\neq\emptyset\}=0.$$
Therefore, we have
$$\limsup_{n\rightarrow \infty} \frac{1}{n} H_{\mu_{\sigma}}(\mathcal{P}^n|\mathcal{Q}^n)\leq \limsup_{n\rightarrow \infty} \frac{1}{n} \sup_{Q_n \in \mathcal{Q}^n} \log \#\Gamma(Q_n).$$
This completes the proof of Theorem  \ref{Thm:conditional-entropy}.
\end{proof}

\subsection{The end of the proof of Theorem~\ref{Thm:entropy-bound-partition}}
Let $\mu$ be an ergodic measure with exactly one positive Lyapunov exponents and $\mathcal{Q}$ satisfies ${\rm Diam}(\mathcal{Q})<\varepsilon$ and $\mu(\partial \cQ)=0$.
For each $\tau>0$, by \eqref{eq: Entropymu}, \eqref{eq:firstpart} and Theorem \ref{Thm:conditional-entropy}, we have 
\begin{align*}
h_\mu(f)\le h_\mu(f, \mathcal{Q})+ \frac{1}{r-1+\alpha}\big(\frac{1}{q}\int \log\|D_xf^q\|{\rm d}\mu(x)-\lambda^+(\mu,f)+\frac{1}{q}\big)+\frac{\log(2 q\Upsilon \cdot C_{r,\alpha})}{q}+\tau.
\end{align*}
The arbitrariness of $\tau>0$ implies the desired result of Theorem \ref{Thm:entropy-bound-partition}.

\section{Proof of Theorem \ref{Thm:B-Invar}}
In this section, we give the proof of Theorem \ref{Thm:B-Invar}. 
Before doing so, we first prove that the sum of the positive Lyapunov exponents is upper semi-continuous.

\begin{Proposition}\label{Pro:usc-sum-positive}
	In any dimension, if $f_n \rightarrow f$ in the $\mathcal C^1$ topology and $\mu_n \rightarrow \mu$ as $n\rightarrow \infty$, where $\mu_n$ is an invariant measure of $f_n$ and $\mu$ is an invariant measure of $f$, then $$\limsup_{n\to\infty}\lambda_\Sigma^+(\mu_n,f_n)\le\lambda_{\Sigma}^+(\mu,f).$$
\end{Proposition}
\begin{proof}
	Assume that $\dim M=d$.
	We give a formula for $\lambda_\Sigma^+(\mu,f)$
	$$\lambda_\Sigma^{+}(\mu,f):=\lim_{n\rightarrow \infty} \frac{1}{n}\int \max_{1\leq k\leq d}\log^{+} \|\wedge^k D_xf^n\| {\rm d}\mu(x).$$
	It is clear that $\varphi_n(x):=\max_{1\leq k\leq d}\log^{+} \|\wedge^k D_xf^n\|$ is continuous.
	\begin{Claim}
		$\{\varphi_n\}$ is a sequence of sub-additive functions.
	\end{Claim}
	\begin{proof}[Proof of the Claim]
		Given $x\in M$ and $n,m\in \NN$, there is $1\le k\le d$ such that
		\begin{align*}
			\varphi_{n+m}(x)&=\log^+\|\wedge^k D_xf^{n + m}\|\leqslant\log^+\|\wedge^k D_{f^{m}(x)}f^{n}\|+\log^+\|\wedge^k D_xf^{m}\|\\
			&\leqslant\varphi_{n}(f^{m}(x))+\varphi_{m}(x).
		\end{align*}
		This concludes the claim.
	\end{proof}
	By Kingman's sub-additive ergodic theorem, one has that
	\begin{equation}\label{eq:sub}
		\lambda_\Sigma^{+}(\mu,f):=\inf_{n\ge 1}\frac{1}{n}\int \max_{1\leq k}~ \log^{+} \|\wedge^k D_xf^n\|{\rm d}\mu(x).
	\end{equation}
	For each $n>0$, the map $(f,\mu)\mapsto \frac{1}{n}\int \max_{k}~\log^{+} \|\wedge^k Df^n\|{\rm d}\mu$ is continuous.
	Then, the upper semi-continuity of $(\mu,f)\mapsto \lambda_{\Sigma}^{+}(\mu,f)$ follows from the above formula.
\end{proof}

\begin{proof}[Proof of Theorem \ref{Thm:B-Invar}]
	Assume that $\dim M=3$.
	Fix $q\in \NN$, $\Upsilon>0$, $r\in \NN$ and $\alpha\in [0,1]$ with $r+\alpha>1$.
	Let $f:M\to M$ be a $\cC^{r,\alpha}$ diffeomorphism with $\|f\|_{\cC^{r,\alpha}}<\Upsilon$.
	Let $\mu$ be an $f$-invariant measure. 
	We decompose $\mu$ into three mutually singular parts:
	$\mu=\alpha^{+}\mu^{+}+\alpha^{-}\mu^{-}+\alpha^0\mu^{0}$, where $\mu^{+}$ has exactly one positive Lyapunov exponent; $\mu^{-}$ has exactly one negative Lyapunov exponent; and $\mu^{0}$ is the remaining part, and these measures are both $f$-invariant measures.
	
    Let $\varepsilon_q$ be as in Theorem \ref{Thm:entropy-bound-partition}, \ref{Thm:entropy-bound-partition-inverse}, and let $\mathcal Q$ be a finite partition with ${\rm Diam}(\mathcal Q)<\varepsilon$.
	Note that almost every ergodic component of $\mu^{+}$ has exactly one positive Lyapunov exponent, while  almost every ergodic component of $\mu^{-}$ has exactly one negative Lyapunov exponent.
	Then, by Theorem \ref{Thm:entropy-bound-partition}, Theorem \ref{Thm:entropy-bound-partition-inverse} and by the ergodic decomposition,
	\begin{align*}
		h_{\mu^{+}}(f)\leq h_{\mu^{+}}(f,\cQ)+ \frac{1}{r-1+\alpha}\big(\frac{1}{q}\int \log\|Df^q\|{\rm d}\mu^{+}-\lambda^+(\mu^{+},f)+\frac{1}{q}\big)+\frac{\log(2 q\Upsilon \cdot C_{r,\alpha})}{q} \\
		h_{\mu^{-}}(f)\leq h_{\mu^{-}}(f,\cQ)+ \frac{1}{r-1+\alpha}\big(\frac{1}{q}\int \log\|Df^{-q}\|{\rm d}\mu^{-}+\lambda^-(\mu^{-},f)+\frac{1}{q}\big)+\frac{\log(2 q\Upsilon \cdot C_{r,\alpha})}{q}. 
	\end{align*}
	Note that $\lambda^+_{\Sigma}(\mu^{+},f)=\lambda^+(\mu^{+},f)$ and
	\begin{align*}
		\frac{1}{q}\int \log\|D_xf^q\|{\rm d}\mu^{+}(x)-\lambda^+(\mu^{+},f) \leq \frac{1}{q}\int \max_{1\leq k\leq 3} \log^{+}~\|\wedge^k D_xf^q\|{\rm d}\mu^{+}(x)-\lambda^+_{\Sigma}(\mu^{+},f).
	\end{align*}
	By $\mu^{\pm}$ are $f$-invariant, and the sum of all Lyapunov exponents is equal to the integral of the log of the Jacobian, we have 
	\begin{align*}
		  &\frac{1}{q}\int \log\|D_xf^{-q}\|{\rm d}\mu^{-}(x)+\lambda^-(\mu^{-},f) \\
		=~&\frac{1}{q}\int \log\|D_{f^q(x)}f^{-q}\|{\rm d}\mu^{-}(x)+\lambda^-(\mu^{-},f) \\
		=~&\frac{1}{q}\int -\big(\log |\Jac D_{x}f^q|- \log \|\wedge^2 D_xf^q\|\big) {\rm d}\mu^{-}(x)+\lambda^-(\mu^{-},f) \\
		=~&\frac{1}{q}\int \log \|\wedge^2 D_xf^q\| {\rm d}\mu^{-}(x)-\lambda^{+}_{\Sigma}(\mu^{-},f) \\
		\leq ~& \frac{1}{q}\int \max_{1\leq k\leq 3}\log^{+}~ \|\wedge^k D_xf^q\| {\rm d}\mu^{-}(x)-\lambda^{+}_{\Sigma}(\mu^{-},f).
	\end{align*}
	By Ruelle’s inequality \cite{Rue78}, the metric entropy of $\mu^0$ is zero, and by \eqref{eq:sub}
	$$\frac{1}{q}\int \max_{1\leq k\leq 3}\log^{+}~ \|\wedge^k D_xf^q\| {\rm d}\mu^{0}(x)-\lambda^{+}_{\Sigma}(\mu^{0},f)\geq 0.$$
	Therefore, by the affine property of metric entropy and partition entropy we have
	\begin{align*}
		h_{\mu}(f) \leq  h_{\mu}(f,\cQ)+\frac{1}{r-1+\alpha}\big(\frac{1}{q}\int \max_{1\leq k\leq 3}\log^{+}~ \|\wedge^k Df^q\| {\rm d}\mu-\lambda^{+}_{\Sigma}(\mu,f)+\frac{1}{q}\big)+\frac{\log(2 q\Upsilon C_{r,\alpha})}{q}. 
	\end{align*}
	This completes the proof of Theorem \ref{Thm:B-Invar} in the case $\dim M=3$.
	The case $\dim M=2$ is completely analogous.
\end{proof}

\section{Proof of Theorem~\ref{Thm:main-perturbation-USC}}
We now turn to the proof of Theorem \ref{Thm:main-perturbation-USC}.
The proof of Theorem \ref{Thm:main-perturbation-USC} is essentially a direct application of Theorem \ref{Thm:B-Invar}.
\begin{proof}[Proof of Theorem \ref{Thm:main-perturbation-USC}]
	Under the assumptions of Theorem \ref{Thm:main-perturbation-USC}, let $\{f_n\}$ be a sequence of diffeomorphisms converging to $f$ in the $\cC^{1,\alpha}$ topology. 
	Choose $\Upsilon>0$ such that $\|f_n\|_{\cC^{1,\alpha}}<\Upsilon$ for every $n\geq 0$.
	For convenience, for each $q\in \NN$ we denote 
	$$q(g,\nu):=\frac{1}{q}\int \max_{1\leq k\leq 3}~\log^{+} ~ \|\wedge^k D_xg^q\|~{\rm d}\nu(x),\quad ~g\in \{f,~f_1,\cdots\},~\nu \in \{\mu,~\mu_1,\cdots\}.$$
	Note that $q(f,\mu) \to \lambda_{\Sigma}^{+}(f,\mu)$ as $q\to \infty$.
	Given $\varepsilon>0$, we first choose $q>0$ such that
	$$\frac{1}{q\cdot \alpha} + \frac{\log (2q\Upsilon C_{\alpha}) }{q}\leq \varepsilon~~\text{and}~~~q(f,\mu)-\lambda_{\Sigma}^{+}(f,\mu)<\alpha\varepsilon,$$
	where $C_{\alpha}$ is the constant given by Theorem \ref{Thm:B-Invar}.
	Let $\varepsilon_q>0$ be the finite constant given by Theorem \ref{Thm:B-Invar}, and let $\cQ$ be a finite partition whose boundary is piecewise smooth, such that $\diam \cQ<\varepsilon_q$ and $\mu(\partial \cQ)=0$.
	Then, by Theorem \ref{Thm:B-Invar}, for each $n>0$ we have 
	$$h_{\mu_n}(f_n)\leq h_{\mu_n}(f_n,\cQ)+\frac{1}{\alpha}\big(\frac{1}{q}\int \max_{1\leq k\leq 3}~\log^{+} ~ \|\wedge^k D_xf_n^q\|~{\rm d}\mu_n(x)-\lambda^{+}_{\Sigma}(\mu_n,f_n)\big)+\varepsilon.$$
	Letting $n\to \infty$, we have $q(f_n,\mu_n)\to q(f,\mu)$ and, since $\mu(\partial \cQ)=0$
	$$\limsup_{n \rightarrow\infty} h_{\mu_n}(f_n,\cQ) \leq  h_{\mu}(f,\cQ).$$
	Therefore, by Theorem \ref{Thm:B-Invar} and  the choice of $q$
	\begin{align*}
	\limsup_{n \rightarrow\infty} h_{\mu_n}(f_n)-h_{\mu}(f)&\leq \limsup_{n \rightarrow\infty} h_{\mu_n}(f_n)-h_{\mu}(f,\cQ) \\
	&\leq \alpha^{-1}\left(q(f,\mu)-\liminf_{n \to\infty} \lambda^{+}_{\Sigma}(\mu_n,f_n)\right) +\varepsilon \\
	&\leq \alpha^{-1}\left(\limsup_{n\to \infty} \lambda^{+}_{\Sigma}(\mu,f)- \lambda^{+}_{\Sigma}(\mu_n,f_n)\right) +2\varepsilon.
	\end{align*}
	This completes the proof of Theorem \ref{Thm:main-perturbation-USC}.
\end{proof}

\section{Proof of Corollary \ref{Cor:SPR-3}} \label{SEC:PofCoE}
We now assume that $M$ is a three-dimensional compact manifold and $f:M\rightarrow M$ is a $\mathcal{C}^{1}$ diffeomorphism.
For $x\in M$, by the sub-additive ergodic theorem, the Lyapunov exponents without multiplicity
$$\lambda^{+}(x,f)\geq \lambda^{c}(x,f) \geq \lambda^{-}(x,f)$$
 are well-defined on a set with full measure for all invariant measures, where
 \begin{align*}
 	&\lambda^{+}(x,f):=\lim_{n \rightarrow \infty} \frac{1}{n} \log^{+} \|D_xf^n\|=\max\{0, \lambda_{1}(x,f)\} \\
 	&\lambda^{-}(x,f):=\lim_{n \rightarrow \infty} -\frac{1}{n} \log^{+} \|D_xf^{-n}\|=-\lambda^{+}(x,f^{-1}).
 \end{align*}
and 
$$\lambda^{c}(x,f):=\lim_{n \rightarrow \infty} \frac{1}{n} \log {\rm Jac}(D_xf^n)-\lambda^{+}(x,f)-\lambda^{-}(x,f)=-\lambda^{c}(x,f^{-1}).$$
For an invariant measure $\mu$, we define
$$\lambda^{+}(\mu,f)=\int \lambda^{+}(x,f) {\rm d}\mu(x),~~\lambda^{c}(\mu,f)=\int \lambda^{c}(x,f) {\rm d}\mu(x),~~\lambda^{-}(\mu,f)=\int \lambda^{-} (x,f) {\rm d}\mu(x).$$
Note that $\mu \mapsto \lambda^{+}(\mu,f)$ is upper semi-continuous and $\mu \mapsto \lambda^{-}(\mu,f)$ is lower semi-continuous.
We denote $\lambda_{\Sigma}^{-}(\mu,f):=-\lambda_{\Sigma}^{+}(\mu,f^{-1})$.
\begin{Lemma}\label{Lem:ELE}
	Let $\{\mu_n\}$ be a sequence of ergodic measures satisfying $\mu_n \rightarrow \mu$ and $\lambda^{+}_{\Sigma}(\mu_n,f)\rightarrow \lambda^{+}_{\Sigma}(\mu,f)$. Then, we have the following 
	\begin{itemize}
		\item if $\lim\limits_{n\rightarrow \infty} \lambda^{+}(\mu_n,f)=\lim\limits_{n\rightarrow \infty} \lambda^{+}_{\Sigma}(\mu_n,f)$, then $\lambda^{c}(x,f)\leq 0$ for $\mu$-almost every $x$;
		\item if $\lim\limits_{n\rightarrow \infty} \lambda^{-}(\mu_n,f)=\lim\limits_{n\rightarrow \infty} \lambda^{-}_{\Sigma}(\mu_n,f)$, then $\lambda^{c}(x,f)\geq 0$ for $\mu$-almost every $x$.
	\end{itemize}
\end{Lemma}
\begin{proof}
	For the first statement, note that
	$$\lambda^{+}_{\Sigma}(\mu,f)=\lim\limits_{n\rightarrow \infty} \lambda^{+}_{\Sigma}(\mu_n,f)=\lim\limits_{n\rightarrow \infty} \lambda^{+}(\mu_n,f)\leq  \lambda^{+}(\mu,f)\leq \lambda^{+}_{\Sigma}(\mu,f).$$
	Then, we have $ \lambda^{+}(\mu,f)=\lambda^{+}_{\Sigma}(\mu,f)$, which means that 
	$$\int \lambda^{+}(x,f) {\rm d}\mu(x) =\int \lambda^{+}(x,f) +\max\{\lambda^{c}(x,f),0\} {\rm d}\mu(x).$$
	This implies that
	$$\int \max\{\lambda^{c}(x,f),0\} {\rm d}\mu(x)=0,$$
	which in turn yields $\lambda^{c}(x,f)\leq 0$ for $\mu$-almost every $x$.
	For the second statement, since
	$$\lim_{n\rightarrow \infty} \lambda^{+}_{\Sigma}(\mu_n,f)+\lambda^{-}_{\Sigma}(\mu_n,f)=\lim_{n\rightarrow \infty}  \int \log {\rm Jac}(D_xf)~{\rm d}\mu_n(x)=\lambda^{+}_{\Sigma}(\mu,f)+\lambda^{-}_{\Sigma}(\mu,f), $$
	one has $\lim_{n\rightarrow \infty}\limits \lambda^{+}_{\Sigma}(\mu_n,f^{-1})=\lambda^{+}_{\Sigma}(\mu,f^{-1})$. 
	Applying the first part to $f^{-1}$ gives $\lambda^{c}(x,f) \geq 0$ for $\mu$-almost every $x$. This completes the proof.
\end{proof}
\begin{proof}[Proof of Corollary \ref{Cor:SPR-3}]
	By the variational principle, we can choose a sequence of ergodic measures $\{\mu_n\}_{n>0}$ such that $h_{\mu_n}(f)\rightarrow h_{\rm top}(f)$. 
	Passing to a sub-sequence, we may assume $\mu_n \rightarrow \mu$. 
	Then, by Theorem \ref{Thm:main-perturbation-USC}, we conclude that $\mu$ is a measure of maximal entropy.
	
	We prove the second statement by contradiction.
	Suppose there exists a sequence of positive numbers $\{\chi_n\}_{n>0}$ decreasing to zero, i.e., $\chi_n \searrow 0$ as $n\rightarrow \infty$, and for each $n>0$, there exists an ergodic measure $\mu_n$ of maximal entropy such that $\lambda_{i}(\mu_n,f)\in [-\chi_{n},\chi_{n}]$ for some $1 \leq i \leq 3$.
	
	By Ruelle's inequality \cite{Rue78}, we know that $\lambda^{+}(\mu_n,f)>\frac{h_{\rm top}(f)}{2}$ and $\lambda^{-}(\mu_n,f)<\frac{-h_{\rm top}(f)}{2}$. 
	Passing to a sub-sequence, we may assume that either
	$$-\chi_{n}\leq \lambda^{c}(\mu_n,f)\leq 0,~\forall n>0~~\text{or}~~0\leq \lambda^{c}(\mu_n,f) \leq \chi_{n},~\forall n>0,$$
	and $\mu_{n} \rightarrow \mu$ as $n\rightarrow +\infty$.
	Note that $h_{\mu_n}(f) \rightarrow h_{\rm top}(f)$ and $\mu$ is a measure of maximal entropy.   
	
	We now show that $\lambda^{c}(x,f)=0$ for $\mu$-almost every $x$.
	We present the proof for the case where $0 \leq \lambda^{c}(\mu_n,f)\leq \chi_n$ for every $n>0$; the other case follows similarly.
	In this case, we have 
	\begin{align*}
		&\lim_{n \rightarrow\infty} \lambda_{\Sigma}^{+}(\mu_n,f)=\lim_{n \rightarrow\infty} \lambda^{+}(\mu_n,f)+\lambda^{c}(\mu_n,f)=\lim_{n \rightarrow \infty} \lambda^{+}(\mu_n,f);\\
		&\lim_{n \rightarrow \infty} \lambda_{\Sigma}^{-}(\mu_n,f)=\lim_{n \rightarrow \infty} \lambda^{-}(\mu_n,f).
	\end{align*}
	Hence, by Lemma \ref{Lem:ELE} one has $\lambda^{c}(x,f)=0$ for $\mu$-almost every $x$, which contradicts the assumption that $\mu$ is a hyperbolic measure.
	
	We now prove the "Moreover" part.	
	For each sequence of ergodic measures  $\{\mu_n\}_{n>0}$ with $\mu_n\rightarrow \mu$ and $h_{\mu_n}(f)\rightarrow h_{\rm top}(f)$, if there exists two sub-sequence $\{n_j\}_{j>0}$ and $\{n_i\}_{i>0}$ such that 
	\begin{itemize}
		\item $\mu_{n_i}$ has exactly one positive Lyapunov exponents for every $i>0$;
		\item $\mu_{n_j}$ has exactly one negative Lyapunov exponents for every $j>0$.
	\end{itemize}
	Then, by Lemma \ref{Lem:ELE}, we have $\lambda^{c}(x,f)\geq 0$ and $\lambda^{c}(x,f)\leq 0$ for $\mu$-almost every $x$, which contradicts the assumption that $\mu$ is hyperbolic. 
	Therefore,  one can choose $i\in \{1,2\}$ and $N\in \NN$ such that $\lambda_{i}(\mu_n,f)>0>\lambda_{i+1}(\mu_n,f)$ for every $n>N$, and
	$\lambda_{i}(x,f)>0>\lambda_{i+1}(x,f)$  for $\mu$-almost every $x$. 
	
	Since $\mu$ is a measure of maximal entropy, and all Lyapunov exponents of ergodic measures of maximal entropy lie outside the interval $[-\chi,\chi]$.
	We have that $\lambda_{i}(x,f)>\chi>-\chi>\lambda_{i+1}(x,f)$ for $\mu$-almost every $x$. 
\end{proof}

\appendix

\section{The proof of the $\mathcal C^{r,\alpha}$ reparametrization lemma}\label{Sec:alpha-reparametrization}
Now we are going to prove Lemma~\ref{Lem:local-reparametrization}. 
We only consider the case $\alpha\in(0,1]$. 
Otherwise, it is the case stated in Burguet's paper.
The proof is parallel to Burguet's work in \cite{Bur24P}.
\subsection{Lemmas from calculus}
\subsubsection{Higher Order Leibniz Rule}
Consider two $\mathcal C^r$ linear operator-valued functions $Q$ and $R$ defined on a open set $X\subset \mathbb{R}^{k}$, where for each $ $z$, Q(z): \mathbb{R}^{n} \rightarrow \mathbb{R}^{m}$ and $R(z): \mathbb{R}^{m} \rightarrow \mathbb{R}^{p}$ are linear operators. 
Recall the higher order Leibniz Rule:
$$D^{r}_{z}(x \mapsto Q(x) \circ R(x))=\sum_{i = 0}^{r}\tbinom{r}{i}(D^{r - i}_{z}Q)\circ(D^{i}_{z}R).$$
Thus, we obtain the estimate (see \eqref{eq:normd} for the definition of $\|\cdot\|_{0}$ and $\|\cdot\|_{\alpha}$)
$$\|D^{r}(Q\circ R)\|_0\le 2^r\max_{k=0,\cdots,r}\|D^kQ\|_0\max_{k=0,\cdots,r}\|D^kR\|_0.$$
For the $\alpha$-norm, recall that
$$\|D(Q\circ R)\|_\alpha\le \|DQ\|_\alpha\|R\|_0+\|DR\|_\alpha\|Q\|_0.$$
Notice that for any $1\leq k <r$ and $\alpha\in(0,1]$,  $\|D^k R\|_\alpha$ can be bounded by $\|D^{k+1} R\|_0$. 
Thus, given $r\in\mathbb N$ and $\alpha\in(0,1]$, there exists a constant $C_{L,r,\alpha}>0$ such that  
\begin{equation}\label{e.higher-holder}
\|D^{r}(Q\circ R)\|_\alpha\le C_{L,r,\alpha} \cdot \max\{\max_{k=0,\cdots,r}\|D^kQ\|_0,\|D^rQ\|_\alpha\} \cdot \max\{\max_{k=0,\cdots,r}\|D^kR\|_0,\|D^rR\|_\alpha\}.
\end{equation}

\subsubsection{Fa\`a di Bruno's formula}

For the estimate on the $\alpha$-norm, one has the following result:
For $\mathcal C^\alpha$ composable continuous maps $\varphi$ and $\psi$, one has that
$$\|\varphi\circ\psi\|_\alpha\le \|\varphi\|_\alpha\|\psi\|_1^\alpha.$$

For $\mathcal C^r$ functions $F:\mathbb{R}^m\rightarrow\mathbb{R}^p$ and $G:\mathbb{R}^n\rightarrow\mathbb{R}^m$, we define $H = F\circ G:\mathbb{R}^n\rightarrow\mathbb{R}^p$. 
We consider the higher-order derivative of $H$. 
For multi-indexes $\alpha = (\alpha_1,\alpha_2,\cdots,\alpha_n)$ and $\beta = (\beta_1,\beta_2,\cdots,\beta_m)$, where $\alpha_i \in \mathbb{N} \cup \{0\}$ and $\beta_i \in \mathbb{N} \cup \{0\}, 1\leq i \leq n, 1\leq j \leq m$. Define 
$$|\alpha|=\alpha_1+\alpha_2+\cdots+\alpha_n.$$
Recall the Fa\`a di Bruno's formula (the Higher Order Chain Rule)
$$\partial^{\alpha}H(x)=\sum_{\beta:|\beta| = |\alpha|}\frac{\alpha!}{\beta!} \cdot \partial^{\beta}F(G(x)) \cdot \sum_{\substack{\sum_{i = 1}^{m}| \gamma^i |=|\alpha|  \\ |\gamma^i|=|\beta_i |}}\frac{\beta!}{\prod_{i = 1}^{m}\gamma^i!}\prod_{i = 1}^{m}\partial^{\gamma^i}G_i(x)
$$
where 
\begin{itemize}
\item $\gamma^i = (\gamma^i_1,\cdots,\gamma^i_n)$ is a multi-index, $G=(G_1,\cdots,G_m)$ and $G_i: \mathbb{R}^n \rightarrow \RR$ for $1\leq i \leq m$;
\item $\alpha!=\prod_{i=1}^{n}\alpha_i ! $. Similarly for $\beta!$ and $\gamma^i!$.
\end{itemize}

One can rewrite as in \cite[Page 1038]{Bur24}: for any $\beta\in\mathbb N^m$ with $|\beta|\le|\alpha|\le r$, there is a universal polynomial $P_{\beta}\left((\partial^{\gamma}G_{i})_{\gamma,i}\right)$ such that 
$$\partial^{\alpha}(F \circ G)=\sum_{\beta \in \mathbb{N}^n, |\beta|\leq|\alpha|} (\partial^{\beta}F)\circ G\times P_{\beta}\left((\partial^{\gamma}G_{i})_{\gamma,i}\right),$$
with $\gamma\in\mathbb N^n$ and $|\gamma|\le |\alpha|$.
We summarize the key estimate we need from Fa\`a di Bruno's formula.
\begin{Lemma}\label{Lem:Bruno}
Given $n,m,p \in \mathbb{N}$, $r\in \mathbb{N}$ and $\alpha\in(0,1]$, there is a constant $C_{B,r,\alpha}>0$ such that for any $\mathcal C^{r,\alpha}$ function $u: \mathbb{R}^m \rightarrow \mathbb{R}^p$, for any function $v:  \mathbb{R}^n \rightarrow \mathbb{R}^m$ satisfying
$$\max_{i=1,\cdots,r}\|D^j v\|_{0}\le 1,~~~\|D^rv\|_\alpha\le 1$$
then one has that
$$\max_{1\leq s\leq r}\|D^s(u\circ v)\|_{0}\le C_{B,r,\alpha}  \max_{1\leq s\leq r}\|D^su\|_{0}~~{\rm and}~~\|D^r(u\circ v)\|_\alpha \le C_{B,r,\alpha} \max\{\max_{1\leq s\leq r}\|D^su\|_{0},\|D^r u\|_\alpha\}.$$
\end{Lemma}

\subsubsection{The Kolmogorov-Landau's inequality}
We have the following $\mathcal C^{r,\alpha}$ version from \cite[Lemma 6]{Bur12}.

\begin{Lemma}\label{Lem:KL-Holder}
Given $r\in\mathbb N$ and $\alpha\in(0,1]$, there is a constant $C_{K,r,\alpha}>0$ such that for any $\mathcal C^{r,\alpha}$ function $\varphi$, one has that
$$\forall k=0,1,\cdots,r,~~~~~~~\|D^k\varphi\|_0\le C_{K,r,\alpha}(\|\varphi\|_0+\|D^r\varphi\|_\alpha).$$
\end{Lemma}

\subsubsection{Taylor's expansion}\label{Subsec:Taylor}

When one considers a $\mathcal C^{r,\alpha}$ map $\varphi: X \rightarrow \mathbb{R}^m$, where  $X \subset \mathbb{R}^n$ is an open concave set, we use the following Taylor expansion at $x\in X$
$$\varphi(x+a)=\sum_{k = 0}^{r}\frac{1}{k!} [D^k_{x} \varphi] (a)^k + R_{r}(x,a),$$
where $a\in \mathbb{R}^n$ with $x+a\in X$, $(a)^k=(a,\cdots,a)\in (\mathbb{R}^n)^{k}$ and 
$$R_{r}(x,a)=\frac{1}{(r - 1)!}\int_{0}^{1}(1- t)^{r - 1} \big([D^r_{(x+ta)} \varphi-D^r_{x} \varphi](a)^r \big) {\rm d}t.$$
Using the H\"older condition, one has that
$$\|  R_{r}(x,a) \|  \leq  \frac{1}{r!} \cdot \|D^r\varphi\|_\alpha \cdot \| a \|^{\alpha+r}.$$

\subsection{Construction of the reparametrizations}
Consider $\Omega>0$ as in the statement of Lemma~\ref{Lem:local-reparametrization}. 
There is $\varepsilon_\Omega>0$ such that for every $\mathcal C^{r,\alpha}$ diffeomorphism $g$ satisfying 
$$ \|g\|_{\mathcal{C}^{r,\alpha}}:= \max \{\max_{k=1,\cdots,r}\|D^{k}g\|_{0},~\|D^r g\|_{\alpha}\}<\Omega$$
 and every $\varepsilon\in(0,\varepsilon_\Omega)$, one has that
$$\forall x\in M,~~~\|D^s g^x_{2\varepsilon}\|_0\le 3\varepsilon\|D_xg\|,~~ s=1,\cdots,r,~~~\textrm{and}~\|D^r g^x_{2\varepsilon}\|_\alpha\le 3\varepsilon\|D_xg\|.$$
where
$$g^x_{2\varepsilon}:=g\circ \exp_x(2\varepsilon\cdot):~T_x M(1)\to M.$$
For a $\mathcal C^{r,\alpha}$ strongly $\varepsilon$-bounded curve $\sigma:[-1,1]\to M$ and $x\in \sigma_{\ast}:=\sigma([-1,1])$, define
\footnote{In a local chart, $g^x_{2\varepsilon}$ has the following form: $g^x_{2\varepsilon}(v)=g(x+2\varepsilon \cdot v)$, $\sigma^x_{2\varepsilon}$ has the following presentation: $\sigma^x_{2\varepsilon}(t)=\frac{1}{2\varepsilon}(\sigma(t)-x);$
and  $g\circ\sigma(t)=g(x+2\varepsilon \frac{1}{2\varepsilon}(\sigma(t)-x))=g^x_{2\varepsilon}\circ \sigma_{2\varepsilon}^x(t).$}
$$\sigma^x_{2\varepsilon}:=\frac{1}{2\varepsilon}\exp_x^{-1}\circ\sigma : [-1,1]\rightarrow T_x M(1).$$

Let $\sigma:[-1,1] \rightarrow M$ be a $\mathcal C^{r,\alpha}$-strongly $\varepsilon$ bounded curve. 
Note that $g\circ\sigma$ may not be bounded anymore.  
To make it to be bounded, one has to compose some reparametrization. 
For an affine reparametrization $\gamma:[-1,1]\rightarrow [-1,1]$ with contraction $b$ such that 
$$\exists t\in[-1,1]~\text{such that}~\sigma(\gamma(t))=y\in M~\text{and}~\lceil \log \|D_y g\| \rceil= \chi,$$
denoting $\Psi(z)=D_zg_{2\varepsilon}^y$ which is a matrix in a local chart, by Equation~\eqref{e.higher-holder}, one has that
\begin{align*}
\left\|D^{r}(g \circ \sigma \circ \gamma)\right\|_\alpha 
\leq ~&b^{r+\alpha} \left\|D^{r}\left(g_{2\varepsilon}^y \circ \sigma_{2\varepsilon}^y\right)\right\|_{\alpha}\\
\leq ~&b^{r+\alpha} \left\|D^{r-1}\left( \Psi(\sigma_{2\varepsilon}^y(t)) \circ D\sigma_{2\varepsilon}^y(t)\right)\right\|_\alpha\\
\leq ~&C_{L,r,\alpha} \cdot b^{r + \alpha} \cdot \|\Psi \circ \sigma_{2\varepsilon}^y\|_{\mathcal{C}^{r-1,\alpha}} \cdot \|D\sigma_{2\varepsilon}^y\|_{\mathcal{C}^{r-1,\alpha}}, 
\end{align*}
where 
$$ \|\Psi \circ \sigma_{2\varepsilon}^y\|_{\mathcal{C}^{r-1,\alpha}}:= \max\{\max_{k=0,\cdots,r-1}\|D^k(\Psi\circ \sigma_{2\varepsilon}^y)\|_0,~\|D^{r-1}(\Psi\circ \sigma_{2\varepsilon}^y)\|_\alpha\};$$
and 
$$ \|D\sigma_{2\varepsilon}^y\|_{\mathcal{C}^{r-1,\alpha}}:= \max\{\max_{k=0,\cdots,r-1}\|D^k(D\sigma_{2\varepsilon}^y)\|_0,~\|D^{r-1}(D\sigma_{2\varepsilon}^y)\|_\alpha\}.$$
By the fact that $\sigma$ is a $\mathcal C^{r,\alpha}$ curve and strongly $\varepsilon$-bounded, one has that 
$$\forall 1\leq k \leq r,~~\|D^k\sigma_{2\varepsilon}^y\|_0\le \frac{1}{2\varepsilon} \|D^k\sigma\|_{0}\leq \frac{1}{\varepsilon} \|D\sigma\|_{0}\leq  1,~~\|D^r\sigma_{2\varepsilon}^y\|_\alpha\le \|D^r\sigma\|_{\alpha}\leq  \frac{1}{\varepsilon} \|D\sigma\|_{0}\leq 1.$$
Thus, by Lemma~\ref{Lem:Bruno},
 there exists $C_{B,r,\alpha}>0$ such that
$$\max_{k = 0,1,\cdots,r-1}\|D^k (\Psi\circ \sigma_{2\varepsilon}^y)\|_0 \le C_{B,r,\alpha} \cdot \max_{k = 1,\cdots,r}\|D^k g_{2\varepsilon}^y\|_{0} \leq  3C_{B,r,\alpha} \varepsilon \cdot \|D_yg\|.$$
and 
$$\|D^{r-1}(\Psi\circ \sigma_{2\varepsilon}^y)\|_\alpha \le  3C_{B,r,\alpha}\cdot \varepsilon \cdot \|D_yg\|.$$
Thus, we have
\begin{align*}
\| D^{r}(g \circ \sigma \circ \gamma)\|_\alpha 
&\leq 3 C_{B,r,\alpha} \cdot b^{r+\alpha} \| D_yg\| \cdot \| D\sigma\|_{0}\\
&\leq 3 C_{B,r,\alpha} \cdot b^{r-1+\alpha}{\rm e}^{\chi^{+}}\| D(\sigma\circ\gamma)\|_{0}\\
&\leq {\rm e}^{\chi - 10}\| \cdot D(\sigma\circ\gamma)\|_{0}
\end{align*}
where we take $b=(3C_{B,r,\alpha}{\rm e}^{\chi^{+}-\chi + 10})^{\frac{-1}{r-1+\alpha}}$.
Thus, if we want to cover the interval, we need at most $b^{-1}+1$ such affine maps.

We use the H\"older form of Taylor's expansion as in Subsection~\ref{Subsec:Taylor}. For $D(g\circ\sigma\circ\gamma)$, we consider the Taylor polynomial $P$ at $0$ of degree $r-1$, we have the estimate:
$$\|P-D(g\circ\sigma\circ\gamma)\|_0\le {\rm e}^{\chi - 10}\| D(\sigma\circ\gamma)\|_{0}$$

For convenience, we define
$$I(\chi,\chi^{+}):=\{t\in[-1,1]:\lceil \log \| D_{\sigma(t)}g \| \rceil=\chi^+,~\lceil\log\|D_{\sigma(t)}g|_{T_{\sigma(t)} \sigma_*}\|\rceil=\chi\},$$
For every $s\in I(\chi,\chi^{+})$, we have $ {\rm e}^{\chi-1}\|D(\sigma\circ\gamma)(s)\|<\|D(g\circ\sigma\circ\gamma)(s)\|\le {\rm e}^\chi\|D(\sigma\circ\gamma)(s)\|$.
Since $\gamma$ is affine and $\sigma$ is strongly $\varepsilon$-bounded,  it follows that $ {\rm e}^{-1}\|D(\sigma\circ\gamma)\|_0\le \|D(\sigma\circ\gamma)(s)\|\le \|D(\sigma\circ\gamma)\|_0$. 
Thus, one has
$$ {\rm e}^{\chi-2}\|D(\sigma\circ\gamma)\|_0<\|D(g\circ\sigma\circ\gamma)(s)\|\le {\rm e}^\chi\|D(\sigma\circ\gamma)\|_0$$
Therefore, we have that 
 \begin{align*}
\| P(s) \| &\leq \| D(g \circ \sigma \circ \gamma)(s) \| + {\rm e}^{\chi - 10} \| D(\sigma \circ \gamma) \|_0 \\
&\leq  {\rm e}^{\chi} \| D(\sigma \circ \gamma) \|_0 + {\rm e}^{\chi - 10} \| D(\sigma \circ \gamma) \|_0 \\
&\leq  {\rm e}^{\chi} \| D(\sigma \circ \gamma) \|_0 (1 +  {\rm e}^{-10}) \\
&\leq  {\rm e}^{3}  {\rm e}^{\chi} \| D(\sigma \circ \gamma) \|_0
\end{align*}
and
\begin{align*}
\| P(s) \| &\geq \| D(g \circ \sigma \circ \gamma)(s) \| -  {\rm e}^{\chi - 10} \| D(\sigma \circ \gamma) \|_0 \\
&\geq  {\rm e}^{\chi - 2} \| D(\sigma \circ \gamma) \|_0 -  {\rm e}^{\chi - 10} \| D(\sigma \circ \gamma) \|_0 \\
&\geq  {\rm e}^{\chi - 2} \| D(\sigma \circ \gamma) \|_0 (1 -  {\rm e}^{-8}) \\
&\geq  {\rm e}^{-3}  {\rm e}^{\chi} \| D(\sigma \circ \gamma) \|_0.
\end{align*}
By the Bezout theorem, there is a constant $C_{B,r}$ depending only on $r$ such that the semi-algebraic set 
$\left\{s\in[-1,1]:~\|P(s)\|\in  ({\rm e}^{-3}  {\rm e}^{\chi} \| D(\sigma \circ \gamma) \|_0, {\rm e}^{3}  {\rm e}^{\chi} \| D(\sigma \circ \gamma) \|_0) \right\}$
is the disjoint union of closed intervals $\{J_i\}_{i\in I}$ with $\# I\le C_{B,r}$.
Moreover, for each $t\in J_i$, one has that
\begin{align*}
\|D(g \circ \sigma \circ \gamma)(t)\| &\leq \|P(t)\| + {\rm e}^{\chi - 10} \|D(\sigma \circ \gamma)\|_0 \\
&\leq {\rm e}^3 {\rm e}^{\chi} \|D(\sigma \circ \varphi)\|_0 + {\rm e}^{\chi - 10} \|D(\sigma \circ \gamma)\|_0 \\
&\leq {\rm e}^4 {\rm e}^{\chi} \|D(\sigma \circ \varphi)\|_0.
\end{align*}
 Let $\gamma_i:[-1,1]\rightarrow J_i$ be the composition of $\gamma$ with an affine map from $[-1,1]$ to $J_i$. 
 
 Now we apply the Kolmogorov-Landau inequality (Lemma~\ref{Lem:KL-Holder}), one has that there is a constant $C_{K,r,\alpha}$ such that for any $1\le s\le r$,
\begin{align*}
\left\|D^s(g\circ\sigma\circ\gamma_i)\right\|_0 &\leq C_{K,r,\alpha} \left(\left\|D^r(g\circ\sigma\circ\gamma_i)\right\|_{\alpha}+\left\|D(g\circ\sigma\circ\gamma_i)\right\|_0\right)\\
&\leq C_{K,r,\alpha} \frac{|J_i|}{2} \left(\left\|D^r(g\circ\sigma\circ\gamma)\right\|_{\alpha}+\sup_{t\in J_i} \left\|D(g\circ\sigma\circ\gamma)(t)\right\|\right)\\
&\leq  C_{K,r,\alpha} \frac{|J_i|}{2} \big({\rm e}^{\chi - 10} \left\|D(\sigma\circ\gamma)\right\|_0+{\rm e}^4{\rm e}^\chi \left\|D(\sigma\circ\gamma)\right\|_0\big)\\
&\leq C_{K,r,\alpha} \frac{|J_i|}{2} {\rm e}^5  {\rm e}^\chi \left\|D(\sigma\circ\gamma)\right\|_0.
\end{align*}
Cut each $J_i$ into $[(1000{\rm e}^5 C_{K,r,\alpha})^{2/\alpha}]+1$ interval $\widetilde J_i$ with the same length. 
Let $\widetilde\gamma_i$ the composition of $\gamma$ with an affine map from $[-1,1]$ onto $\widetilde J_i$. 
By the construction, the number of affine maps of the form $\widetilde \gamma_i$ such that the union of their images can cover $I(\chi^{+},\chi)$ is at most
$$(3C_{B,r,\alpha}{\rm e}^{\chi^{+}-\chi + 10})^{\frac{1}{r-1+\alpha}}\cdot C_{B,r}  \cdot(1000{\rm e}^5 C_{K,r,\alpha})^{2/\alpha}:=C_{r,\alpha} \cdot \exp \left(\frac{\chi^{+}-\chi}{r-1+\alpha} \right).$$
We are now going to check that $g\circ\sigma\circ\widetilde\gamma_i$ is bounded.
We first consider $2\le s\le r$:
\begin{align*}
\left\|D^s(g\circ\sigma\circ\widetilde{\gamma}_i)\right\|_0 &\leq (1000{\rm e}^5C_{r,\alpha,K})^{-2} \left\|D^s(g\circ\sigma\circ\gamma_i)\right\|_0\\
&\leq \frac{1}{6}(1000{\rm e}^5 C_{K,r,\alpha} )^{-1} \frac{|J_i|}{2}{\rm e}^5{\rm e}^\chi \left\|D(\sigma\circ\gamma)\right\|_0\\
&\leq \frac{1}{6}(1000C_{K,r,\alpha} )^{-1} \frac{|J_i|}{2} \min_{s\in J_i} \left\|D(g\circ\sigma\circ\gamma)(s)\right\|\\
&\leq \frac{1}{6}(1000C_{K,r,\alpha} )^{-1} \frac{|J_i|}{2} \min_{s\in \widetilde{J}_i} \left\|D(g\circ\sigma\circ\gamma)(s)\right\|\\
&\leq \frac{1}{6} \left\|D(g\circ\sigma\circ\widetilde{\gamma}_i)\right\|_0.
\end{align*}
Now we check for $r+\alpha$. By using the fact that $\| D^{r}(g \circ \sigma \circ \gamma)\|_\alpha\le {\rm e}^{\chi - 10}\| D(\sigma\circ\gamma)\|_{0}$, it is very similar to the above estimate:
\begin{align*}
\left\|D^r(g\circ\sigma\circ\widetilde{\gamma}_i)\right\|_{\alpha} &\leq (1000{\rm e}^5 C_{K,r,\alpha})^{-2} \left\|D^r(g\circ\sigma\circ\gamma_i)\right\|_{\alpha}\\
&\leq \frac{1}{6}(1000{\rm e}^5C_{K,r,\alpha})^{-1} \frac{|J_i|}{2}{\rm e}^\chi \left\|D(\sigma\circ\gamma)\right\|_0\\
&\leq \frac{1}{6}(1000C_{K,r,\alpha})^{-1} \frac{|J_i|}{2} \min_{s\in J_i} \left\|D(g\circ\sigma\circ\gamma)(s)\right\|\\
&\leq \frac{1}{6}(1000C_{K,r,\alpha})^{-1} \frac{|J_i|}{2} \min_{s\in \widetilde{J}_i} \left\|D(g\circ\sigma\circ\gamma)(s)\right\|\\
&\leq \frac{1}{6} \left\|Dg\circ\sigma\circ\widetilde{\gamma}_i\right\|_0.
\end{align*}
This completes the proof of Lemma \ref{Lem:local-reparametrization}.

\subsection*{Data availibility statement}
No data availability statement is required, as no experimental data is involved.

\subsection*{Ethics declarations}
\noindent \textbf{Conflict of interest} \
The authors state that there is no conflict of interest.

\section*{Acknowledgements}
The authors would like to thank professor J. Buzzi and professor D. Burguet for their comments.

\vskip 5pt

\flushleft{\bf Chiyi Luo} \\
\small School of Mathematics and Statistics, Jiangxi Normal University, Nanchang,   330022, P. R. China\\
\textit{E-mail:} \texttt{luochiyi98@gmail.com}\\
\flushleft{\bf Dawei Yang} \\
\small School of Mathematical Sciences,  Soochow University, Suzhou, 215006, P.R. China\\
\textit{E-mail:} \texttt{yangdw@suda.edu.cn}\\
\end{document}